\documentclass[letterpaper,twoside,onecolumn]{ametsocV6.1}

\title{Analysis of tidal flows through the Strait of Gibraltar using Dynamic Mode Decomposition}

\authors{Sathsara Dias \aff{a},
Sudam Surasinghe\aff{a,b,f},\correspondingauthor{Sudam Surasinghe, Sudam.surasinghe@yale.edu}
K.G.D. Sulalitha Priyankara\aff{a,c}, 
Marko Budi\v{s}i\'{c}\aff{a}, 
Larry Pratt\aff{d},
Jos\'{e} C.~Sanchez-Garrido\aff{e},
and Erik M. Bollt\aff{f}
}

\affiliation{\aff{a}{Department of Mathematics, Clarkson University, 8 Clarkson Ave,Potsdam, NY 13699, USA}\\
\aff{b}{Department of Ecology \& Evolutionary Biology , Yale University, 165 Prospect Street, New Haven, CT 06511, USA }\\
\aff{c}{Department of Mathematics, Alabama A\&M University, 4900 Meridian St N, Huntsville, AL, 35811, USA }\\
\aff{d}{Woods Hole Oceanographic Institution, 266 Woods Hole Road MS\#21, Woods Hole, MA, 02543-1535, USA}\\
\aff{e}{Department of Applied Physics II, University of Málaga, Avda. Cervantes 2, Málaga, 29071, Spain}\\
\aff{f}{Department of Electrical and Computer Engineering, Clarkson University, 8 Clarkson Ave, Potsdam, NY 13699, USA }
}
 
\abstract{The Strait of Gibraltar is a region characterized by intricate oceanic sub-mesoscale features, influenced by topography, tidal forces, instabilities, and nonlinear hydraulic processes, all governed by the nonlinear equations of fluid motion. In this study, we aim to uncover the underlying physics of these phenomena within 3D MIT general circulation model simulations, including waves, eddies, and gyres. To achieve this, we employ Dynamic Mode Decomposition (DMD) to break down simulation snapshots into Koopman modes, with distinct exponential growth/decay rates and oscillation frequencies. Our objectives encompass evaluating DMD's efficacy in capturing known features, unveiling new elements, ranking modes, and exploring order reduction. We also introduce modifications to enhance DMD's robustness, numerical accuracy, and robustness of eigenvalues. DMD analysis yields a comprehensive understanding of flow patterns, internal wave formation, and the dynamics of the Strait of Gibraltar, its meandering behaviors, and the formation of a secondary gyre, notably the Western Alboran Gyre, as well as the propagation of Kelvin and coastal-trapped waves along the African coast.
In doing so, it significantly advances our comprehension of intricate oceanographic phenomena and underscores the immense utility of DMD as an analytical tool for such complex datasets, suggesting that DMD could serve as a valuable addition to the toolkit of oceanographers.}

\begin{document}

\maketitle

%
%
%
\statement
A challenge arising in all branches of geophysics is making sense of increasingly large data sets describing  complex processes.  In physical oceanography this data usually comes from high-resolution models and observations.  Reduced order models such as the dynamic mode decomposition (DMD) have the potential to identify key processes and interactions through the synthesis of a data-based model with relatively few of degrees of freedom.  Because of its connection with Koopman operator theory, DMD is also able, in principle, to deal with data describing nonlinear processes.   In this work, we test the ability of DMD to describe the essential physics exhibited in a complicated data set produced by a model of the ocean circulation in Strait of Gibraltar and Western Mediterranean, a region that contains striking, time-dependent, and often nonlinear features, mostly driven or modulated by tides. We carefully describe the connection to Koopman theory and the step-by-step DMD algorithm, and we present a procedure for singling out the most robust of the DMD modes.   We then show how particular modes contain information about specific physical processes including modulation of the two-layer exchange flow in the strait, the generation of internal waves, meandering of the Atlantic Jet, and the generation of coastal-trapped waves.  Some of these features are already well-known, but others such as the jet meandering have received little attention.     
%
%
%


\section{Introduction}

The interpretation of geophysical data often involves a search for dominant features or physical processes.
Examples include spatial distributions dominated by a small number of empirical orthogonal functions, vertical structure dominated by the lowest few vertical dynamical modes, processes responsible for peaks in energy spectra, or multiscale phenomena describable by simple scaling laws.
Once such a feature is identified, it may be possible to gain deeper insight through examination using reduced order models, including idealized models grounded in the governing equations of motion or models that are constructed entirely from field data or data generated by large models of general circulation and climate.

The Dynamic Mode Decomposition (DMD)~\citep{schmid2010,rowley2009jfm} and its various iterations have, over the past decade, been shown to identify dominant, interpretable features of a data set that lead to low-order, linear models describing the time evolution of the features.
Among other competing approaches, such as Proper Orthogonal Decomposition (POD)~\citep{smith2005NlinDyn}, Empirical Orthogonal Functions (EOF)~\citep{lorenz1956empirical,weare1982} and others (see~\citep{taira2020,taira2017} for  modern comparative reviews), DMD has been distinguished by an interpretation that casts it as an approximation of the more general framework employing the linear Koopman operator for analysis of nonlinear dynamical systems~\citep{mezic2004,mezic2005,korda2018, korda2018a,mezic2013,mezic2019,surasinghe2021randomized} or, simply, Koopman analysis.

Conceptually, Koopman analysis augments the set of variables on which the model depends by a set of variables computed by \emph{all} functions on the original state space (\emph{observables)}.
As a consequence, any nonlinear term present in the original equations is given its own "synthetic" variable; as this leads to an explosion of variables, the effect is that the linearity of the model is obtained in return for the dimension of the model growing.
At first glance, transformation to a linear system of infinite order would appear to be the opposite of order reduction, but hope is that the infinite set can be represented to a good approximation by a large, but finite, dynamical system and that the solutions will exhibit low-order behavior.
The numerical algorithms, like those in the DMD family, strive to compute eigenvalues and eigenvectors of this linear, infinite-dimensional linear system without explicitly trying to reconstruct its system matrix.
The computed eigenvalues indicate the character of time variation of individual components: their interpretation is the same as in basic stability theory of linear systems.
The eigenvectors approximate Koopman modes whose interpretation depends on the choice of observed variables; when input variables are measurements of velocities in the computational domain, as is common in fluid mechanics, Koopman modes correspond to spatial profiles that grow, decay, and/or oscillate according to a single complex frequency (Koopman eigenvalue).

Since the seminal publication that first connected Schmid's DMD with Koopman analysis~\citep{rowley2009jfm}, many iterations of the DMD algorithm have been developed with standard goals common to many numerical algorithms such as improved robustness, improved efficiency, reduction of bias and noise, and improved convergence properties.
Several reviews and edited volumes address various aspects of the interaction between Koopman analysis, numerics, and applications~\citep{mezic2013,brunton2022,budisic2012chaos,mauroy2020}.
While the proliferation of algorithms and applications lead to significant improvements on the original technique, it made it more challenging for those looking to test or adapt these algorithm to find an entry point.

Applications of DMD to fluid mechanics are arguably best studied, as early papers focused primarily on demonstrating the DMD on simulated fluid systems, typically in aerodynamic flows.
However, DMD has not been widely applied to large oceanographic data sets, either those including models or observations.

The central goal of our work is to determine how well DMD works when applied to a region of the ocean that exhibits a variety of time-dependent, 3D physical features including waves, eddies, and gyres, some of which are generated by strongly nonlinear processes.
How well do the Koopman modes capture the features that we already know about, and does the method reveal features that we do not know about?
It is of additional interest to rank the numerous modes that are computed and determine what reduction in order can be achieved.
Here we will discuss several alternatives.

The oceanic region selected includes the Strait of Gibraltar and the Western Alboran Sea.
The data consists of tide-resolving MITgcm output over a 6-day period and is detailed in Section~\ref{sec:model}. Prominent known features include the Western Alboran Gyre, the Atlantic Jet, the hydraulically controlled exchange flow in the Strait of Gibraltar, the time-dependent hydraulic jump that occurs to the west of the Camarinal Sill, and the internal waves that radiate eastward into the Mediterranean and are generated by the collapse of the hydraulic jump.
The collapse is a strongly non-linear process caused by the flooding of the semi-diurnal tide into the Mediterranean Sea.
Therefore, the standard harmonic analysis~\citep{pawlowicz2002,foreman1989} that consists of fitting the data to a predetermined set of steady oscillation frequencies is not likely to capture the full breadth of dynamical behavior.

As a secondary goal, we aim to demonstrate how to apply the DMD to a data set in a straightforward and easy-to-implement fashion, while including several modern modifications aimed at improving numerical performance of the algorithm.
We demonstrate how certain modeling choices affect the outcomes of the algorithm, whether by interpreting changes in individual DMD modes, or by evaluating the robustness and correctness of the reduced-order models in bulk.

\section{Model and analyzed data}\label{sec:model}

The analysis in this paper utilizes simulated current velocities in the Strait of Gibraltar region using the MIT general circulation model (MITgcm~\citep{Marshall1997}). The numerical calculations were conducted with a non-hydrostatic version of the MITgcm regional application described in \citep{sg2013,sg2015}. For a detailed description of the model the reader is referred to these papers; here we summarize the main features. The model domain covers the Alboran Sea and the Gulf of Cadiz (connected by the Strait of Gibraltar) and was discretized with a curvilinear orthogonal grid of variable resolution. Horizontal resolution varies from 0.5 km in the Strait itself to 3-4 km towards the adjacent basins. In the vertical direction, the model incorporates 46 unevenly distributed z-levels with maximum resolution at the near surface ($\Delta z=5 \si{\meter}$). The circulation was driven at the free surface by wind stress, heat, and freshwater fluxes from NCEP/NCAR reanalysis~\citep{Kalnay96}. Initial and lateral open boundary values for temperature, salinity, and velocity were derived from CMEMS reanalysis products~\citep{sotillo15}. Tides were introduced by laterally forcing the model with barotropic velocities associated to the main semidiurnal (M2, S2, N2, K2) and diurnal (K1, O1, P1, Q1) tidal constituents. The outputs analyzed here are a subset of 144 hourly snapshots of 3-D current velocity –vector components and magnitude–, covering a total of 6 days. The size of each individual snapshot is 190x96x32. A summary of the descriptions and symbols employed for velocity-related quantities throughout the article can be found in Table~\ref{tab:Notations}.

\begin{table*}[hbt]
\begin{center}
\begin{tabular}{lc}
\topline
Description & Symbol \\
\midline
Zonal velocity (Longitude-wise velocity component) & $U_x$ \\
Meridional velocity (Latitude-wise velocity component) & $U_y$ \\
Velocity along depth & $U_z$ \\
Local horizontal speed & $U_s=\sqrt{U_x^2+U_y^2}$ \\
Time average of $U_x$ & $\overline{U_x}$ \\
Time average of $U_y$ & $\overline{U_y}$ \\
Time average of $U_z$ & $\overline{U_z}$ \\
Time average of $U_s$ & $\overline{U_s}$ \\
Magnitude/speed of the mean (time averaged) velocity vector & $S_{\overline{U}}= \sqrt{\overline{U_x}^2+\overline{U_y}^2}$ \\
\botline
\end{tabular}
\end{center}
\caption{Summary of descriptions and symbols for velocity-related quantities as used throughout the article.}\label{tab:Notations}
\end{table*}

\begin{figure*}[htb]
  \centering
  \includegraphics[width=\textwidth]{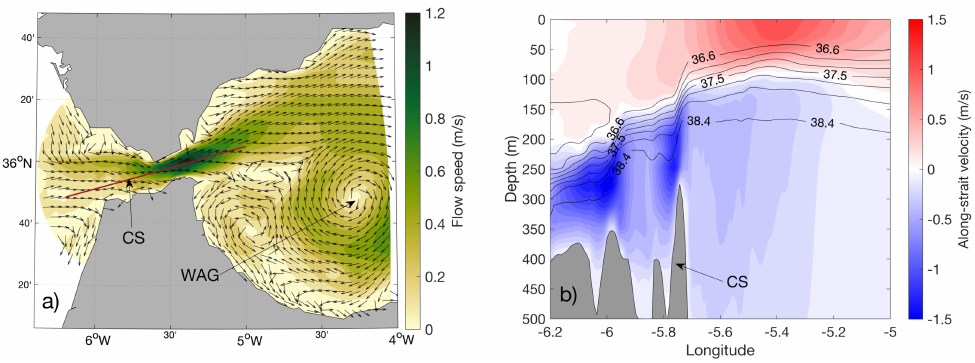}
  \caption{a) Time average surface velocity. Colors indicate the flow speed ($S_{\overline{U}}$); arrows indicate the flow direction. b) Time average along-strait velocity in a vertical cross section following the Strait’s axis (red line in panel a). The location of the Camarinal Sill (CS) is labeled in both panels.}
  \label{fig:1_R1}
\end{figure*}

The main features of the flow are illustrated here from the model outputs and described based on the analysis of \citet{sanchez2011numerical,sanchez-garrido2013a}. The time average velocities show a strong surface current entering the Alboran Sea through the Strait of Gibraltar exceeding the speed of 1 m/s in the region of maximum narrowing –Tarifa Narrows– (Figure~\ref{fig:1_R1}a). This current, referred to as the Atlantic jet, enters the Alboran Sea in a northeastward direction and later follows a surrounding path around a the anticyclonic Western Alboran Gyre (WAG). A vertical cross section along the strait axis reveals the two-layer structure of the exchange flow, with a countercurrent of Mediterranean water flowing towards the Atlantic beneath the Atlantic jet (Figure~\ref{fig:1_R1}b). The flow exhibits plunging isopycnals (or equivalently here isohalines) and accelerated Mediterranean currents at the western flank of the Camarinal sill (CS; the main sill of the strait). This supercritical overflow is preceded by a hydraulic control section at the CS crest. Further west, in the so-called Tangier basin (\(\sim\ang{6.0}\text{W}-\ang{5.8}\text{W}\)), the flow returns to subcritical through an internal hydraulic jump. Here large dissipation of kinetic energy and mixing occurs, the latter suggested by the widening experienced by the isohalines at \(\sim \ang{5.9}\text{W}\). Similar internal hydraulic transitions take place further west at Espartel Sill (\(\sim \ang{6}\text{W}\)). 

\begin{figure*}[htb]
  \centering
  \includegraphics[width=\textwidth]{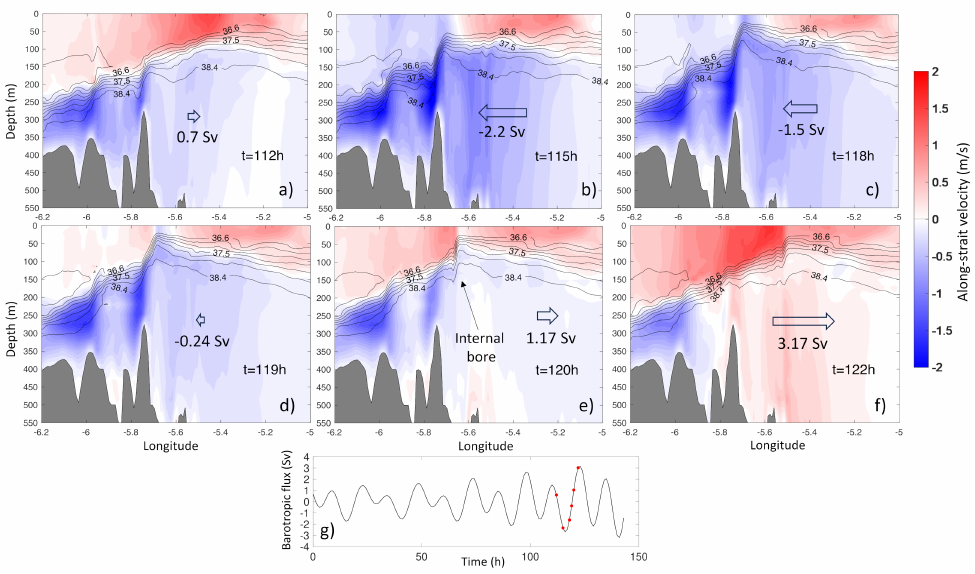}
  \caption{Snapshots of along-strait velocity (colors) and salinity (contours) during a tidal cycle (a-f). Arrows in each panel indicate the strength and direction of the barotropic flux. The time series of barotropic flux is shown in the bottom panel (g), with dots corresponding to the frames a-f.}
  \label{fig:2_R2}
\end{figure*}

Both Atlantic and Mediterranean layers are strongly modulated by tides. The volume transport of the mean exchange is ~0.85 Sv, while barotropic tidal flows can reach peak values of 5 Sv, with 3 Sv associated to the M2 tidal constituent alone (\cite{lafuente2000tide}). Tidal currents are strong enough to interrupt either the inflow or outflow during a certain phase of the tidal cycle, especially at CS. This and other time-dependent effects are illustrated in Figure~\ref{fig:2_R2}. During the flood tide (westward tidal flow) the Mediterranean outflow strengthens, blocking the Atlantic flow over the CS and at the western part of the strait (Figure~\ref{fig:2_R2}b-c). Other effects of the westward tidal flow are the intensification of the hydraulic jump located at the lee side of CS –and Espartel sill, further west– and the rising of the isohalines as one approaches CS from the east (Figure~\ref{fig:2_R2}c). The westward shoaling of the interface between Atlantic and Mediterranean waters to the east of the sill (roughly 5.6°W to 5.0°W) reduces the phase velocity of one of the long internal waves supported by the flow. Eventually, a disturbance that attempts to propagate towards the Mediterranean Sea becomes stationary and a new critical section, called approach control, emerges east of the sill. In the plots shown in Figure~\ref{fig:2_R2}c-d, such a control is located around 5.6°W. Note that there is no strong visual feature that marks the position of the approach control but a sign that a control is present is a small hydraulic jump lying between it and the critical section at CS.  There is a hint of such a feature in Figure~\ref{fig:2_R2}c, where the isopynals in the upper layer rapidly rise immediately westward of 5.6°W. More details on the two-dimensional structure of the sill flow and quantitative evidence for the approach control are shown in the supplementary material.

\begin{figure*}[htb]
  \centering
  \includegraphics[width=\textwidth]{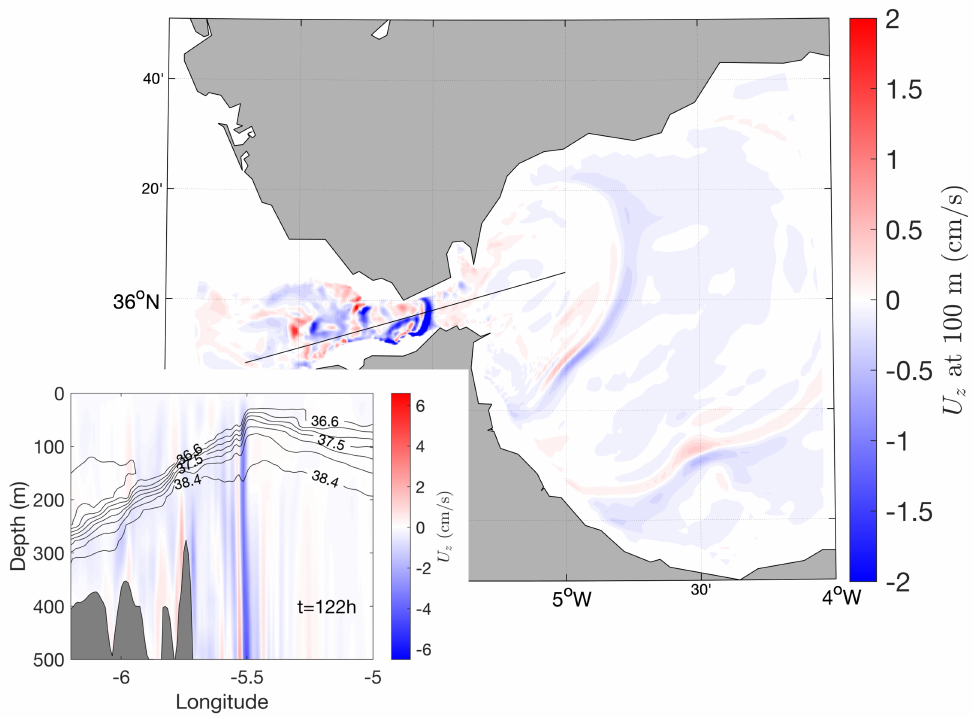}
  \caption{Snapshot of vertical velocity at 100 m depth. The two wave patterns in the Alboran Sea correspond to two internal bores generated during consecutive semi-diurnal tidal cycles. A third trailing wave front in the strait narrows corresponds to a nascent internal bore progressing eastward. Inset: cross section of vertical velocity and salinity (contours), showing the vertical structure of the bore at 5.5ºW.}
  \label{fig:R3}
\end{figure*}

\begin{figure*}[htb]
  \centering
  \includegraphics[width=\textwidth]{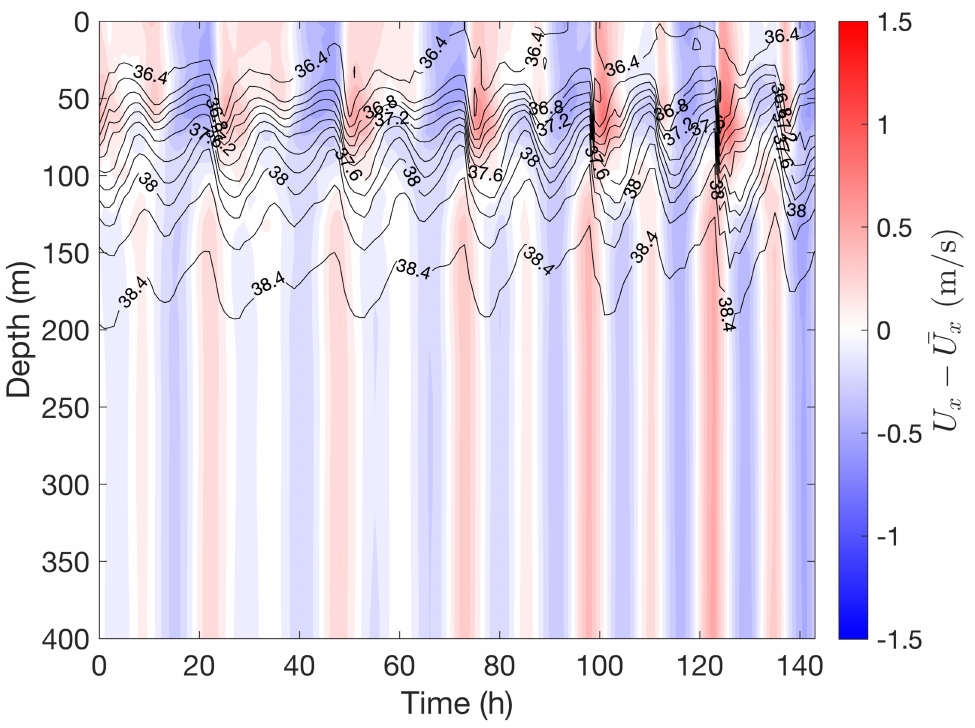}
  \caption{Time-depth plot of zonal velocity anomaly (anomaly with respect to the time mean) at 5.4º W, 36.0º N (eastern side of the Strait of Gibraltar).}
  \label{fig:R4}
\end{figure*}

The approach control is short-lived; it is lost a couple of hours after it was formed as the westward flow slackens, allowing for the release of an upstream internal hydraulic jump formed between the CS and the approach control. The moving internal hydraulic jump –or internal bore– propagates eastward at a speed of 1-2 m/s (Figure~\ref{fig:2_R2}e-f), breaking down into a group of solitary waves (\cite{sanchez-garrido2008}). The details of the disintegration process require horizontal resolution ~100 m, which is finer than the 500 m resolution of the model, so although a bore is produced it is not captured in detail. The internal bore has undular characteristics and can be seen in the vertical displacements of the isohalines and from the model vertical velocities, showing strong downwelling at the wavefront and upwelling at the rear. (Figure~\ref{fig:R3}). The generation of the internal bore is linked to the mentioned approach control, which is established only under sufficiently strong westward tidal currents. The condition for the generation of the bore is fulfilled in almost every semidiurnal tidal cycle except during neap tides, when the process becomes diurnal due to the diurnal inequality of the semidiurnal tidal currents ($t<96$ h in the present case; Figure~\ref{fig:2_R2}g). Lastly, it should be noted that the tidal dynamics of the Strait lead also to internal tides. Like the bore, internal tides originate at CS and radiate eastward, but they have tidal periodicity and are longer than the bore. A time-depth plot of salinity and horizontal velocity anomaly ($U_x-\overline{U_x}$) at the eastern part of the Strait reveals that internal tides are primarily semidiurnal and have the vertical structure of a first-mode internal disturbance (Figure~\ref{fig:R4}).  

\section{Dynamic Mode Decomposition of the simulated flow}

\subsection{Representing dynamics using the Koopman operator}\label{Sec:BasicTheory}

We start by summarizing the basic concepts of the ``Koopmanism'' approach to dynamical systems.
To properly justify this approach, we have to disentangle the concept of ``states of the system'' and ``observations of the system'' which are often conflated in analysis of dynamical models, with notable exceptions of control theory and data assimilation literature.
We present the basic concepts here without much theoretical justification, and point to~\citep{mezic2013,mezic2005,budisic2012chaos} for more detail.

Consider a dynamical system over states \(z \in \mathcal{Z}\)
\begin{align}
  \frac{dz}{dt} &= g(z)\label{eq:ode}
\end{align}

In the context of our applications, states \(z \in \mathcal{Z}\) correspond to all quantities necessary for evolution of the model equations.
In general, they can be thought of as spectral coefficients if model equations are advanced by spectral solvers, grid-based evaluations of quantities in method-of-lines, etc.
Almost always we assume that dimension of \(\mathcal{Z}\) is large, even infinite.
As the time evolution is discretized into steps of duration \(\Delta t\), states \(z(n \Delta t)\) being evolved by a time-invariant, nonlinear solution operator \(\gls{evol}_{0}^{\Delta t} : \mathcal{Z} \to \mathcal{Z}\):
\begin{align}
 z[n] := z(n \Delta t) = \gls{evol} _{0}^{\Delta t}( z[n-1] ).
\label{eq:evolution}
\end{align}
Since \eqref{eq:ode} was time invariant and as \(\Delta t\) remains fixed, we will write simply \(\gls{evol}  :=  \gls{evol} _{0}^{\Delta t}\equiv \gls{evol} _{t}^{t+\Delta t}\).

Access to \(z(t)\) is mediated by a vector of \(D\) observation functions \(x : \mathcal{Z} \to \mathbb{R}^{D}\);
therefore, we only have access to \emph{snapshots}
\begin{equation}
  \label{eq:snapshot}
  \mat{X}[n] := x( z[n] ) \quad \text{ of } \quad \mat{X}(t) := x(z(t)).
\end{equation}
Ocean general circulation models typically yield observation functions \(x(\cdot)\) that are evaluations of velocity components, the salinity, pressure, temperature, and other quantities at spatial grid points.
We do not access \(x(\cdot)\) as functions, but only as a time-series \(\mat{X}[n] = x( z[n] )\) generated by evaluating \(x(\cdot)\) along the otherwise-inaccessible state evolution \(z[n]\).
This means that at any time step \(n\), vector \(\mat{X}[n]\) contains \(D = \# variables \times \# gridpoints \) values, so that the \emph{snapshot matrix} \(\mat{X}\) representing an evolution over \(N\) timesteps is a \(D \times N\) matrix.

The nature of the state space \(\mathcal{Z}\) is important only to the extent that it allows for a definition of the vector space of functions \(\mathcal{X}\) to which elements of \(x(\cdot)\) belong.
In applied contexts, this space is commonly taken to be the space \(\mathcal{X} = L^{p}(\mathcal{Z},m)\) of \(\mathbb{C}\)-valued functions, in particular with \(p=2\) where measure \(m\) is taken either to be a conserved measure of dynamics, or volume (Lebesgue measure) on \(\mathcal{Z}\).
Notice that this implies that \(\dim \mathcal{X} \geq \dim \mathcal{Z}\).  In practice, a function in $\mathcal{X}$ may be interpreted as a physical scalar  measurement of the process over the space $\mathcal{Z}$, whereas it can be accumulated (``averaged")  across the space by $m$ accordingly.

The \emph{Koopman operator} \(\gls{koop}  : \mathcal{X} \to \mathcal{X}\) is the evolution rule advancing time for measurement functions  \(x \in \mathcal{X}\) by composing them with evolution operator \(\gls{evol} \) defined in (\ref{eq:evolution})
  \begin{align}
  \label{eq:koopman}
\gls{koop} x &= x \circ \gls{evol} ,\\
\shortintertext{or, evaluated pointwise,}
\gls{koop} \underbrace{x(z[n])}_{\mat{X}[n]} &= x( \gls{evol} (z[n]) ) = \underbrace{x(z[n+1])}_{\mat{X}[n+1]}.
\end{align}
Operator \(\gls{koop} \) is analogous to \(\gls{evol} \) for states \(z\).
Such evolution implies that \(\mat{X}[n]\) at any point is determined by \(\mat{X}[0]\):
\begin{equation}
  \label{eq:koopman-evolution}
  \mat{X}[n] = \gls{koop} \mat{X}[n-1] = \gls{koop}^{2}\mat{X}[n-2] = \cdots = \gls{koop}^{n}\mat{X}[0].
\end{equation}

Operator \(\gls{koop} \) is linear, i.e., for any two functions \(x, \hat x\) and scalars \(a,b\),
\begin{equation}
  \label{eq:linearity}
  \gls{koop} (a x(\cdot) + b \hat x(\cdot) ) = a x(\gls{evol} (\cdot)) + b \hat x(\gls{evol} (\cdot)) = a \gls{koop}  x + b \gls{koop}  \hat x.
\end{equation}
Notice that this property holds regardless of \(\gls{evol} \) being typically nonlinear.
Effectively, the nonlinearity of \(\gls{evol} \) on \(\mathcal{Z}\) is ``unraveled'' into additional dimensions of \(\mathcal{X}\).

Scalar \(\gls{kdev} \in \mathbb{C}\) and function \(\gls{kef}:\mathcal{Z}\to \mathbb{C}\) of \(\gls{koop} \) that satisfy
\begin{equation}
  \label{eq:eigenpair}
  \gls{koop} \gls{kef} = \gls{kdev} \gls{kef},
\end{equation}
are called an eigenvalue and eigenfunction, respectively (together eigenpair).
If \(\gls{kef} \in \mathcal{X}\) happens to be one of the measurement functions selected to form an element of \(\mat{X}[n]\), the (\ref{eq:eigenpair}) implies that its evolution \(\gls{kef}{}[n] = \gls{kef}( z[n])\) would be given simply through multiplication by the associated eigenvalue
\begin{equation}
  \label{eq:evolution-of-eigenfunction}
  \gls{kef}{}[n] = \gls{koop}^{n} \gls{kef}{}[0] = \gls{kdev}^{n} \gls{kef}{}[0].
\end{equation}
We can introduce \emph{continuous-time eigenvalues}
\begin{equation}
  \label{eq:continuous-time-eigenvalues}
  \gls{kev} = \gls{kdr} + i \gls{kaf} := \frac{1}{\Delta t} \ln \gls{kdev}, \quad \text{ or } \quad \gls{kdev} = \exp(\gls{kev} \Delta t)
\end{equation}
so that the evolution \(\gls{kdev}^{n} = e^{\gls{kev} n \Delta t}\) has an interpretation in terms of the original continuous time scale \(t\) both qualitatively and quantitatively.

Given that \(\gls{evol} \) is linear, the eigenpairs of \(\gls{evol} \) directly relate to the eigenpairs of the Koopman operator \(\gls{koop}\). If \( \lambda \) and \(\gls{kef} \) are an eigenvalue and eigenfunction of \(\gls{evol} \), then \(\gls{koop}\) has the eigenvalue \( e^{\lambda \Delta t} \) with the same eigenfunction \(\gls{kef} \). This relationship effectively connects the discrete-time evolution under \(\gls{evol} \) with the continuous-time spectral characteristics represented by  \(\gls{koop}\).

To aid interpretation, we will for eigenvalues often note their period of oscillation \(\gls{per}\)
\begin{equation}
  \label{eq:period-of-oscillation}
  \gls{per} = 2\pi/\gls{kaf}
\end{equation}
and halving/doubling time \(\gls{dou}\)
\begin{equation}
  \label{eq:doubling-time}
  e^{\gls{kdr}  \gls{dou}} = 2, \quad \gls{dou} := \frac{\ln 2}{\gls{kdr}}
\end{equation}
which are expressed units of time.
Interpretation as doubling vs.~halving time is detected qualitatively from the sign of \(\gls{kdr}\) or \(\gls{dou}\), positive for growing eigenfunctions and doubling time, negative for decaying eigenfunctions and halving time.

By linearity, if any measurement function \(x_{p}\) of the snapshot is a linear combination of eigenfunctions \(\gls{kef}_{k}\), \(x_{p}(\cdot) = \sum_{k=0}^{K} \gls{kmd}_{pk} \gls{kef}_{k}(\cdot)\) with coefficients \(\gls{kmd}_{pk} \in\mathbb{C}\),  \(K \leq \infty\), then its time-trace would be given by
\begin{align}
  \label{eq:evolution-of-the-span}
  \mat{X}_{p}[n] = \gls{koop}^{n} \left(\sum_{k=0}^{K} \gls{kmd}_{pk} \gls{kef}_{k}\right) = \sum_{k=0}^{K} \gls{kmd}_{pk} \gls{kdev}_{k}^{n}\underbrace{\gls{kef}_{k}[0]}_{=: \gls{kcoeff}_{k}}
\end{align}
Expanding \emph{all} elements of vector-valued observation \(x(\cdot)\) into the same set of eigenfunctions, gives rise to
\begin{equation}
  \begin{bmatrix}
| \\ \mat{X}[n] \\ |
  \end{bmatrix} =
  \sum_{k=0}^{K}
  \begin{bmatrix} | \\ \gls{kmd}_{k} \\ | \end{bmatrix}
  \gls{kdev}_{k}^{n}
  \gls{kcoeff}_{k}
\quad  \text{or} \quad
  \begin{bmatrix}
| \\ \mat{X}(t) \\ |
  \end{bmatrix} =
  \sum_{k=0}^{K}
  \begin{bmatrix} | \\ \gls{kmd}_{k} \\ | \end{bmatrix}
  e^{t \gls{kev}_{k} }
  \gls{kcoeff}_{k}
\label{eq:koopman-mode-decomposition}
\end{equation}
where the two formulas are, respectively, continuous- and discrete-time versions of the decomposition, linked at timesteps \(t = n\Delta t\).
The column-vectors \(\gls{kmd}_{k}\equiv \gls{kmd}_{\ast k} \in \mathbb{C}^{D}\) are called \emph{Koopman modes}, and these modes represent the spatial patterns that evolve over time according to the dynamics captured by the eigenvalues.

We note here the similarity of the ``separation-of-variables'' in Koopman mode decomposition (KMD) \eqref{eq:koopman-mode-decomposition} with the evolution employed by Empirical Orthogonal Functions (EOF)~\citep{lorenz1956empirical,Monahan2009}, or Proper Orthogonal Decomposition (POD)~\citep{smith2005NlinDyn}:
\begin{equation}
  \label{eq:pod-eof}
  \mat{X}[n] =   \underbrace{\sum_{k=0}^{K}
  \gls{kmd}_{k}   \gls{kcoeff}_{k}
  \gls{kdev}_{k}^{n}}_{\text{KMD} }
 = \underbrace{\sum_{k=0}^{K} \mat{U}_{k} \sigma_{k} \mat{V}_{k}[n]}_{\text{EOF/POD}}.
\end{equation}

the matrix form
\begin{figure*}[t] 
\begin{align}  \label{eq:koopman-mode-decomposition-matrix}
\begin{bmatrix}
  | &  & | &  \\
  \mat{X}[0] & \cdots & \mat{X}[n] & \cdots \\
  | &  & | &  \\
  \end{bmatrix} &=
  \underbrace{ \begin{bmatrix}
& | & \\
\cdots & \gls{kmd}_{k} & \cdots \\
& | & \\
  \end{bmatrix}}_{\gls{kmd}}
\underbrace{  \begin{bsmallmatrix}
\ddots &  & \\
 & \gls{kcoeff}_{k} &  \\
&  & \ddots\\
  \end{bsmallmatrix} }_{\diag \gls{kcoeff}}
\underbrace{\begin{bmatrix}
   & &  & \\
  \gls{kdev}_{k}^{0} & \cdots & \gls{kdev}_{k}^{n} & \cdots \\
   &  &  & \\
  \end{bmatrix}}_{\gls{kdev}^{\top}}
\end{align}
\end{figure*}

While EOF/POD modes \(\mat{U}_{k}\) are real-valued and orthonormal, \(\mat{U}_{k}^{\top} \mat{U}_{k'} = \delta_{kk'}\), Koopman modes are typically not orthonormal and appear in complex-conjugate pairs to capture oscillatory components.
On the other hand, there is no \emph{a priori} dynamics that can be assumed about temporal evolution of EOF/POD modes \(\mat{V}_{k}[n]\), while Koopman mode evolution is easy to interpret, as it is always a decaying/growing/periodic evolution \(\gls{kdev}_{k}^{n} = e^{\gls{kdr}_{k} n \Delta t}( \cos(\gls{kaf}_{k} n \Delta t) \pm i  \sin(\gls{kaf}_{k} n \Delta t))\).
Further discussion of the connections is available in~\citep[\S 4.3]{tu2014a}.

The dimension of a mode \(\gls{kmd}_{k}\) is equal to the number of observables used to form \(\mat{X}[n]\) because mode \(\gls{kmd}_{k}\) stores the coefficients with which eigenfunction \(\gls{kef}_{k}\) participates in evolution of snapshots \(\mat{X}[n]\).
Koopman modes correspond to spatial profiles that evolve in time according to a single time constant \(\gls{kdev}_{k} = e^{\Delta t \gls{kev}_{k}}\).

When eigenvalues \(\gls{kdev} = e^{\gls{kev}}\) are real-valued, the associated component of the decomposition \eqref{eq:evolution-of-the-span} evolves as
\begin{equation}
  \label{eq:real-valued-modes}
  \sim \gls{kmd}e^{\gls{kev} t}\gls{kcoeff}.
\end{equation}
On the other hand, oscillating (complex-valued) modes appear in  complex-conjugate pairs \(\gls{kmd}, \gls{kmd}_{k'}^{\ast}\) that evolve according to eigenvalues \(\gls{kdev}, \gls{kdev}^{\ast}\).
The associated real-valued component of the decomposition \eqref{eq:evolution-of-the-span} will be in the form
\begin{equation}
  \label{eq:complex-valued-modes}
  \begin{aligned}
  \sim  \gls{kmd}^{\ast} e^{\gls{kdr}t-i \gls{kaf} t} \gls{kcoeff}^{\ast} + \gls{kmd} e^{\gls{kdr}t + i \gls{kaf} t} \gls{kcoeff} & =
  2 \Re\left( \gls{kmd} e^{\gls{kdr}t+i \gls{kaf} t} \gls{kcoeff} \right) \\
&=  2 \abs{\gls{kmd}}\abs{\gls{kcoeff}}e^{\gls{kdr} t} \cos( \gls{kaf}t + \Theta + \theta)
  \end{aligned}
  \end{equation}
  where the polar form of the mode and of the associated coefficient are
\begin{equation}
\gls{kmd} = \abs{\gls{kmd}}e^{i\Theta}, \quad    \gls{kcoeff} = \abs{\gls{kcoeff}}e^{i\theta}.\label{eq:polar-form-of-modes}
  \end{equation}
  Here \(\abs{\gls{kmd}}, \Theta\) are element-wise moduli and phases of \(\gls{kmd}\).

  In summary, while both elements of the conjugate pair are needed for a \emph{linear} decomposition \eqref{eq:evolution-of-the-span}, the formulation in terms of the magnitude and phase of the associated oscillation \eqref{eq:complex-valued-modes} is typically preferred for interpretation of the associated dynamics.

clarifies that Koopman framework amounts to a formalism akin to solving a PDE by separation of variables, where the columns of the modal matrix \(\gls{kmd}\) are \emph{spatial profiles} constant in time, while rows of the \(\gls{kdev}^{\top}\) matrix are the \emph{temporal time traces} according to which each of the modes evolves.
The diagonal matrix \(\diag \gls{kcoeff}\) plays the role of the coefficients used to fit this spatiotemporal decomposition onto the first snapshot, as for the first \(n=0\) column this reduces to
\begin{align}  \label{eq:koopman-initial-condition}
  \mat{X}[0] &=
     \gls{kmd}
  \begin{bsmallmatrix}
\ddots &  & \\
 & \gls{kcoeff}_{k} &  \\
&  & \ddots\\
  \end{bsmallmatrix}
  \begin{bmatrix}
1 \\ 1 \\ \vdots
  \end{bmatrix} = \gls{kmd} \gls{kcoeff}.
\end{align}
Remembering that \(\gls{kcoeff}\) is derived from  eigenfunctions of the Koopman operator~\eqref{eq:evolution-of-the-span} and the initial condition of the evolution \(z(0)\), we realized that the eigenfunctions and initial conditions are entangled in a vector of coefficients \(\gls{kef}{}[0] = \gls{kef}(z(0))\) formed by evaluations of all eigenfunctions at the initial condition of the evolution.
Theoretically, dimension of \(\gls{kef}\) is infinite, although rank of the input data matrix \(\mat{X}\) limits the number of eigenfunctions that are practically accessible.

\subsection{Dynamic Mode Decomposition as numerical Koopman analysis}\label{sec:dmd}
Dynamical Mode Decomposition (DMD)~\citep{schmid2010,rowley2009jfm} is the algorithm that is able to approximate the decomposition~\eqref{eq:evolution-of-the-span} in a non-parametric manner, i.e., without \emph{a priori} specifying the location of eigenvalues.
While the derivation of this algorithm is not new~\citep{tu2014a}, we repeat it here because it clarifies how the choice of parameters affects the outcome, and allows us to comment on modifications we performed in Section~\ref{sec:modifications}.

Starting with the snapshot matrix \(\mat{X} = \mat{X}[n]\), \(n=0, 1, \dots, N-1\), subdivide it into a ``past'' and ``future'':
\begin{figure*}[t] 
\begin{equation}
  \label{eq:x-subdivision}
  \begin{aligned}
\mat{X}_{1}&=\begin{bmatrix} \mat{X}[0] & \mat{X}[1] & \dots & \mat{X}[N-2] & \phantom{\mat{X}[n-1]} \end{bmatrix}\in \mathbb{R}^{m\times (N-1)} \\
\mat{X}_{2}&=\begin{bmatrix} \phantom{\mat{X}[0]} & \mat{X}[1] & \dots & \mat{X}[N-2] & \mat{X}[N-1] \end{bmatrix}\in \mathbb{R}^{m\times (N-1)} \\
\end{aligned}
\end{equation}
\end{figure*}


While evolution \(\mat{X}(z[n]) \mapsto \mat{X}(z[n+1]) = \gls{koop} \mat{X}( z[n] ) \) is connected by the infinite-dimensional Koopman operator, we approximate it using left-multiplication by the DMD matrix \(\gls{dmdk}\)
\begin{align}
  \mat{X}[n+1] &= \gls{dmdk} \mat{X}[n], \\
  \shortintertext{or}
  \mat{X}_{2} &= \gls{dmdk} \mat{X}_{1}.
\label{eq:dmd-dynamics}
\end{align}
This amounts to the assumption that \(\mat{X}[n]\) were generated as iterations of a finite-dimensional linear map. The goal of DMD algorithms is to approximate eigenpairs of \gls{dmdk} without approximating the (dense) matrix \(\gls{koop}\) itself.

Further assuming that this map has a complete set of eigenvectors \(\gls{kmd}_m\) with eigenvalues \(\gls{kdev}_{m} \in \mathbb{C}\), which is the generic case under small random perturbations of matrix elements, we write
 \begin{equation}\label{eq:dmd-solution}
\mat{X}[n] = \sum_{m=1}^{M} \gls{kmd}_m \gls{kcoeff}_{m} \gls{kdev}_{m}^{n},
\end{equation}
 where coefficients \(\gls{kcoeff}_{m}\) are determined by the initial condition
\begin{equation}
  \label{eq:initial-condition}
  \mat{X}[0] = \sum_{m=1}^{M}\gls{kmd}_{m} \gls{kcoeff}_{m} = \gls{kmd}\gls{kcoeff}.
\end{equation}
To make the choice of coefficients unique, we normalize each \(\gls{kmd}_{m}\) to \(\norm{\gls{kmd}_{m}} = 1\).
We refer to \(\gls{kmd}_{m}\) as \emph{DMD modes}, to \(\gls{kdev}_{m}\) as \emph{DMD eigenvalues}, and to \(\gls{kcoeff}_{m}\) as \emph{reconstruction coefficients}.
Magnitudes \(\abs{\gls{kcoeff}_{m}}\) are norms of DMD modes.

Of course, this decomposition of dynamics exactly mirrors \eqref{eq:koopman-mode-decomposition}, so the common thread for all DMD algorithms is computation of eigenvalues and eigenvectors of \gls{dmdk}.
In principle Equation \eqref{eq:dmd-dynamics} can be solved approximately by finding the \(\ell^2\)-optimal solution that minimizes the Frobenius norm
\begin{equation}
\gls{dmdk} := \operatorname*{arg\,min}_{\mat{A}} \norm{\mat{A}\mat{X}_{1} - \mat{X}_{2}}_{F}.\label{eq:l2-dmd}
\end{equation}
When the number of time-steps is smaller than the dimension of a snapshot, making \(\mat{X}\) a ``tall'' matrix, this is equivalent to forming the  Moore--Penrose pseudoinverse \(\mat{X}_{1}^{\pinv}\) and setting
\begin{align}
  \label{eq:dmd-matrix}
  \gls{dmdk} &= \mat{X}_{2}\mat{X}_{1}^{\pinv} \\
  \shortintertext{ which amounts to}
  \gls{dmdk} &= \gls{koop} \mat{U} \mat{U}^{\top},
\end{align}
where \(U\) is the orthonormal basis for the data subspace \(\lspan\{  \mat{X}[n] \}_{n=0}^{N-2}\), so that  \(U U^{\top}\) is the orthogonal projector onto the subspace spanned by first \(N-1\) snapshots.

This route is often impractical.
Storing the square matrix \(\gls{dmdk}\) requires \(m^{2}\) values, even though we are practically interested in some small subset of vectors (DMD modes), each of length \(m\).
Moreover, operations involving this matrix are known to be numerically-fragile ~\citep{laub1985numerical}, so an alternative route is needed.
Finally, solving \eqref{eq:initial-condition} for coefficients \(b\) again requires for a fitting, which can be performed either in the same or different norm/process as fitting of \(\gls{dmdk}\).
Various flavors of the DMD algorithm available in the literature differ in the choices they make to address these questions (among their other differences).

Here we present the simplest version, the so-called ``exact'' DMD~\citep{tu2014a}.
In it, the fitting of \(\gls{dmdk}\) is made in the \(\ell^{2}\)-space of matrices, using the Singular Value Decomposition.
The dimensionality of \(\gls{dmdk}\) is resolved by first projecting dynamics to a lower-dimensional subspace of \(r\) singular vectors \(\mat{U}_{r}\) of $\mat{X}_1$, and then (implicitly) choosing a non-orthogonal basis of DMD modes in that subspace, to recover the requirement that each mode corresponds to a single (complex) time constant~\eqref{eq:continuous-time-eigenvalues}.
Finally, the fitting of \(b\) is again performed in \(\ell^{2}\)-space of vectors.
We point out that the process below does not completely resolve the issue of numerical robustness; for details see~\citep{Drmac2019,Drmac2018}.

First, we compute the truncated singular value decomposition of $\mat{X}_1$
\begin{equation}
  \mat{X}_1 \approx \mat{U}_{r} \mat{\Sigma}_{r} \mat{V}_{r}^{\top}
\end{equation}
where $\mat{U}_{r} \in \mathbb{R}^{m \times r}$, $ \mat{V}_{r} \in \mathbb{R}^{(N-1) \times r}$, $\mat{\Sigma}_{r} \in \mathbb{R}^{r \times r}$ are matrices of first \(r\) left- and right-singular vectors and singular values, respectively.

Equation \eqref{eq:dmd-dynamics} then becomes
\begin{equation}
  \mat{X}_2 \approx \gls{dmdk} \mat{U}_{r} \mat{\Sigma}_{r} \mat{V}_{r}^{\top},
\end{equation}
and therefore,
\begin{equation}
 \gls{romk} :=  \mat{U}_{r}^{\top} \gls{dmdk} \mat{U}_{r}=\mat{U}_{r}^{\top} \mat{X}_2 \mat{V}_{r} \mat{\Sigma}_{r}^{-1}.
\end{equation}
Since \(\mat{U}_{r}\) is a subset of columns of \(U\), this implies that
\begin{equation}
\gls{romk} = \mat{U}_{r}^{\top} \gls{koop} \mat{U}_{r}\label{eq:romk-koopman}
\end{equation}
so that eigenvalues of \(\gls{romk}\) are a subset of eigenvalues of \(\gls{dmdk}\) and \(\gls{koop}\).

The Eigendecomposition
\begin{equation}
 \gls{romk}  \gls{rommod}_k = \gls{kdev}_k \gls{rommod}
\end{equation}
gives a good approximation of Koopman eigenvalues, but \(\gls{rommod}\) are still only \(r\)-dimensional.

Merely ``undoing'' the compression by \(\mat{U}_{r}\)  and setting \(\gls{kmd} \approx \gls{kmd}'_{k} := \mat{U}_{r} \gls{rommod}\) produces eigenvectors that have the correct number of elements, but are constrained to be in the artificially-reduced span \(\mat{U}_{r}\):
\begin{equation}
\mat{U}_{r}\mat{U}_{r}^{\top} \gls{dmdk} \underbrace{\mat{U}_{r} \gls{rommod}_k}_{\gls{kmd}'_{k}} = \gls{kdev}_k \underbrace{\mat{U}_{r}\gls{rommod}}_{\gls{kmd}'_{k}}.
\end{equation}
Matrix $\mat{U}_{r}\mat{U}_{r}^{\top} \gls{dmdk}$ is the DMD matrix whose output is projected onto the lower-dimensional subspace.
Such calculation may be favorable in cases where the size of the input data set prohibits computing more than  \(r\) largest singular vectors~\citep{Stoll2012}, but if the full SVD decomposition is available, \emph{exact} DMD can do better.
From SVD  \(\mat{X}_{1} = \mat{U}\mat{\Sigma} \mat{V}^{\top}\), where number of columns of \(\mat{U}\) and \(\mat{V}\) is equal to \(\operatorname{rank} \mat{X}_{1}\),  the exact DMD prescription makes use of it by setting
\begin{equation}
   \gls{kmd}_k := \mat{X}_2 \mat{V} \mat{\Sigma}^{-1}\gls{rommod}_{k}.
 \end{equation}

 To calculate combination vector \gls{kcoeff}, one solves the initial condition equation \eqref{eq:initial-condition} by \(\ell^{2}\)-minimization
 \begin{equation}\label{eq:b-coefficients}
\mat{X}[0] = \gls{kmd} \gls{kcoeff} \rightarrow \gls{kcoeff} = \gls{kmd}^{\pinv} \mat{X}[0],
\end{equation}
which then ensures that the data \(\mat{X}[0]\) is well-approximated (in the \(\ell^{2}\) sense)
 \begin{equation}\label{eq:dmd-approximation}
\mat{X}[n] \approx \sum_{k}\gls{kmd}_{k} \gls{kdev}^{n}_{k} \gls{kcoeff}_{k} = \sum_{k}\gls{kmd}_{k} e^{\gls{kev}_{k} n \Delta t } \gls{kcoeff}_{k}.
\end{equation}
The full process is summarized in Algorithm~\ref{alg:exact-DMD}.

DMD eigenvalues and modes are not canonically ordered by importance; this is in contrast to POD/EOF modes which are typically ordered by the magnitude of the associated singular value.
Therefore, ordering and selection of DMD modes is an additional step that is typically peformed by taking account of the purpose that the DMD computation serves in the overall modeling process.
Since the mode selection is less standard than mode computation, as it may be specific to the model analyzed, we address it in a separate section.

\begin{algorithm}[htb]
  \caption{Exact DMD algorithm.
For details see Section~\ref{sec:dmd}.}
  \SetAlgoLined \DontPrintSemicolon
  \SetKwInOut{Input}{input}\SetKwInOut{Output}{output}
  \ResetInOut{output}
  \Input{Data matrix \(\mat{X}\), with columns \(\mat{X}[k], n=0,\dots,N-1\)}
  \Input{Snapshot sampling rate \(\Delta t\)}
  \Input{Dimension of the subspace \(r\)}
  \BlankLine
  \(\mat{X}_{1} \gets \mat{X}[0],\dots,\mat{X}[N-2],\quad \mat{X}_{2} \gets \mat{X}[1],\dots,\mat{X}[N-1] \)\;
  \(U, \mat{\Sigma}, V \gets \textrm{svd}(\mat{X}_{1})\) \tcp*{Singular Value Decomposition}
  \(\mat{U}_{r} \gets \mat{U}[1:r,:], \mat{\Sigma}_{r} \gets \mat{\Sigma}[1:r,1:r], \mat{V}_{r} \gets \mat{V}[1:r,:]\) \tcp*{Reduce to \(r\)-dim.\ subspace}
  \(\gls{romk} \gets \mat{U}_{r}^{\top} \mat{X}_2 \mat{V}_{r} \mat{\Sigma}_{r}^{-1}\) \tcp*{DMD matrix}
  \(\gls{kdev}_k, \gls{rommod}_{k} \gets \textrm{eig}(\gls{romk})\) \tcp*{Eigendecomposition}
  \(\gls{kmd}_k \gets \mat{X}_2 V \mat{\Sigma}^{-1}\gls{rommod}_{k}\) \tcp*{Koopman modes}
  \(\gls{kmd}_{k} \gets \gls{kmd}_{k}/\norm{\gls{kmd}_{k}}_{2}\) \tcp*{Mode normalization}
  \(\gls{kev}_{k} \gets (\ln \gls{kdev}_{k})/\Delta t\) \tcp*{Continuous-time DMD eigenvalues}
  \(\gls{kcoeff} \gets \gls{kmd}^{\pinv} \mat{X}[0]\) \tcp*{\(\ell^{2}\) fit of the 1st snapshot}
  \Output{\(\gls{kmd}_k, \gls{kev}_k, \gls{kcoeff}_{k}\)}
  \label{alg:exact-DMD}
\end{algorithm}

\subsection{Improving interpretability, stability, and robustness}\label{sec:modifications}

For the purposes of this paper, we report on DMD analysis of the three-dimensional velocity fields.
This suggests that each snapshot in \(\mat{X}\) contains three velocity components \([u,v,w]\) evaluated on the spatial grid.
We augment this variable vector by the horizontal speed and use
\begin{equation}
  \frac{\sqrt{2}}{2}
  \begin{bmatrix}
U_x \\ U_y \\  \sqrt{2} U_z\\ \sqrt{U_x^{2} + U_y^{2}}
  \end{bmatrix}\label{eq:observable-stack}
\end{equation}
as the snapshot vector.
Since the horizontal speed \(U_s\coloneqq\sqrt{U_x^{2} + U_y^{2}}\) is a nonlinear function of \(U_x,U_y\), computing its DMD decomposition from knowing only decomposition of \(U_x,U_y\) would require combining all available DMD modes, complicating any model reduction.
Including the horizontal speed as an explicit observable avoids doing so.
The \(\sqrt{2}\) factors ensure that the \(\ell^{2}\) norm of any snapshot approximates the integrated \(\ell^{2}\)-norm of the three-dimensional velocity field, avoiding biasing of any component.
Our simulated data has \(144\) snapshots, amounting to a \(6 \si{\day}\) window sampled at \(1\si{\hour}\).
The size of the computational grid is \(190 \times 96 \times 32\); evaluating the four observables on it and stacking as in~\eqref{eq:observable-stack} yields the snapshot matrix of dimensions \(2,334,720 \times 144\).

Several articles have argued for and against removing the mean value of data as the first pre-processing step~\citep{chen2012,hirsh2020,seenivasaharagavan2021}.
If the trajectory is long and the dynamics stationary, in the sense that the spectrum does not change over time, computing the mean
\begin{equation}
  \label{eq:ergodic-mean}
  \gls{kmd}_{0} \coloneqq \frac{1}{N} \sum_{n=1}^{N} X[n]
\end{equation}
and subtracting it from data amounts to orthogonally-projecting the data onto the \(\lambda=1\) eigenfunction of the Koopman operator~\citep{mezic2004}.
This can be thought of as a ``parametric'' estimation of the Koopman mode, as we presume to know \emph{a priori} the eigenvalue \(\gls{kev}=1\) of interest.
However, when the duration of data is short and the dynamics possesses eigenvalues off the unit circle, the practical utility of mean-removal and theoretical implications are less understood.
Recent investigations suggest that mean-removal increases the robustness of the algorithm without compromising the integrity of the analysis~\citep{hirsh2020,seenivasaharagavan2021}.
We have evaluated the outcomes of the algorithm with and without the mean removal and concluded that the approximation quality of non-constant modes is improved in the case where the mean is not removed in pre-processing.
The likely cause is the lack of orthogonality of DMD modes.
Since the mean vector is recurrent with respect to any period \(P\), removing it parametrically as above removes a degree of freedom in subsequent estimation of periodic modes and regression of data on them.
The same problem does not occur in POD/EOF reduced order models, since the constraint of orthogonality of modes prevents the data mean to contribute to any other mode except the ``mean mode''.

As detailed in~\citep{Dawson2016,Hemati2017a}, the baseline DMD algorithm suffers from bias in eigenvalues due to using only \(\mat{X}_{1}\) (past snapshots) to construct the lower-dimensional subspace, amounting to the least-squares (LSQ) regression of the subspace onto the data.
Instead, we use the total-least-squares (TLSQ) correction~\citep{Hemati2017a} by projecting \(\mat{X}_{1},\mat{X}_{2}\) matrices onto a subspace derived from the SVD of a larger \(\mat{Z}\) matrix, formed by vertically stacking \(\mat{X}_{1,2}\) matrices.
Computationally, this step is simple to implement, although the SVD of the matrix \(\mat{Z}\) now requires twice as much memory as before.

We adopt the following additional modifications of the baseline exact DMD based on~\citep{Drmac2018}.
\begin{compactenum}[(a)]
\item Columns of matrices \(\mat{X}_{1}, \mat{X}_{2}\) are normalized by \(\ell^{2}\)-norms of matrix \(\mat{X}_{1}\), improving the condition number and therefore making numerical algebra more robust~\citep[\S 3.4]{Drmac2018}.
\item In our implementation of DMD we replaced MATLAB's stock divide-and-conquer SVD algorithm with the LAPACK QR SVD algorithm which computes the small singular values and the associated singular vectors to a higher accuracy following the recommendation in~\citep[\S 4.1.3]{Drmac2018}.
\item We evaluated the Refined Rayleigh-Ritz Data Driven Modal Decomposition (DDMD RRR) modification to DMD~\citep[\S 3.5]{Drmac2018}, which optimizes the residual of eigenvector equation for \(\gls{romk}\) for the DMD eigenvalues and eigenvectors.
  The changes compared to the baseline exact DMD algorithm were minimal; to avoid a substantial computational cost we did not use the RRR modification here.
\end{compactenum}

The baseline DMD algorithm Algorithm~\ref{alg:exact-DMD} computes \(\gls{kcoeff}\) by minimizing the fit of the DMD model at the initial condition over all complex-valued vectors \(v\), \(\gls{kcoeff} = \argmin_{v}\norm{\gls{kmd}v - \mat{X}[0]}^{2}_{2}\).
Since eigenvalues are already known at that point in the algorithm, in principle any other snapshot, \(n' \not = 0\), could have been chosen for the optimization target as in
\begin{equation*}
  \argmin_{v}\norm*{
 \gls{kmd} \diag\begin{bmatrix} \gls{kdev}_{1}^{n'}  & \gls{kdev}_{2}^{n'}& \dots \end{bmatrix} v - \mat{X}[n']
  }^{2}_{2}.
\end{equation*}
Practically, the size of coefficients \(\gls{kcoeff}\) will change depending on the choice of the time step.
Since all modes were normalized to the unit \(\ell^{2}\) norm, any decaying modes will be assigned relatively higher coefficients when \(n'=0\) is chosen while the growing modes will be assigned relatively higher coefficients when \(n' \gg 0\).
For steady-state data, we should expect that the relative importance of modes does not change, therefore we calculate \(\gls{kcoeff}\) that minimizes the total error across all snapshots

\begin{equation}
\label{eq:computing-b}
\gls{kcoeff} = \argmin_{v} \sum_{n=0}^{N-1}\norm*{
\gls{kmd} \diag\begin{bmatrix} \gls{kdev}_{1}^{n} & \gls{kdev}_{2}^{n} & \dots \end{bmatrix} v - \mat{X}[n]
}^{2}_{2}.
\end{equation}
This equation can be converted to an algebraic linear system through the use of a tensor product of eigenmode and eigenvalue matrices~\citep{dias2024identifying}.
Solving \eqref{eq:computing-b} for all \(n = 0,1,\dots,N-1\) on a computer can require a substantial amount of working memory so, instead,
we use only a small subset of snapshots spread over the entire time domain to compute the sum to ensure minimal biasing of growing/decaying modes, while keeping computation cost at bay.

\begin{algorithm}[htb]
  \caption{Modified exact DMD algorithm.}
  \SetAlgoLined \DontPrintSemicolon
  \SetKwInOut{Input}{input}\SetKwInOut{Output}{output}
  \ResetInOut{output}
  \Input{Data matrix \(\mat{X}\), with columns \(\mat{X}[k], n=0,\dots,N-1\)}
  \Input{Snapshot sampling rate \(\Delta t\)}
  \Input{Dimension of the subspace \(r\)}
  \BlankLine
  \(\mat{X}_{1} \gets \mat{X}[0],\dots,\mat{X}[N-2],\quad \mat{X}_{2} \gets \mat{X}[1],\dots,\mat{X}[N-1]\)\;
  \(\mat{X}_{1}[:,k] \gets \mat{X}_{1}[:,k]/\norm{\mat{X}_{1}[:,k]}_{2},\ \mat{X}_{2}[:,k] \gets \mat{X}_{2}[:,k]/\norm{\mat{X}_{1}[:,k]}_{2}\)\tcp*{Normalization}
  \(\mat{Z} \gets [X1; X2]\)\quad  \tcp*{TLSQ start}
  \(\mat{U}_{\mat{Z}}, \mat{\Sigma}_{\mat{Z}}, \mat{V}_{\mat{Z}} = \textrm{QRsvd}(\mat{Z}), \mat{V}_{\mat{Z},r} \gets \mat{V}_{\mat{Z}}[1:r,:]\) \tcp*{using QR-based SVD}
  \(\mat{X}_{1} \gets \mat{X}_{1} \mat{V}_{\mat{Z},r},\ \mat{X}_{2} \gets \mat{X}_{2} \mat{V}_{\mat{Z},r}\)\tcp*{TLSQ projection}
  \(\mat{U}, \mat{\Sigma}, \mat{V} \gets \textrm{QRsvd}(\mat{X}_{1})\) \tcp*{QR-based SVD}
  \(\mat{U}_{r} \gets \mat{U}[1:r,:], \mat{\Sigma}_{r} \gets \mat{\Sigma}[1:r,1:r], \mat{V}_{r} \gets \mat{V}[1:r,:]\) \tcp*{Reduce to \(r\)-dim.\ subspace}
  \(\gls{romk} \gets \mat{U}_{r}^{\top} \mat{X}_2 \mat{V}_{r} \mat{\Sigma}_{r}^{-1}\) \tcp*{DMD matrix}
  \(\gls{kdev}_k, \gls{rommod}_{k} \gets \textrm{eig}(\gls{romk})\) \tcp*{Eigendecomposition}
  \(\gls{kmd}_k \gets \mat{X}_2 \mat{V} \mat{\Sigma}^{-1}\gls{rommod}_{k}\) \tcp*{Koopman modes}
  \(\gls{kmd}_{k} \gets \gls{kmd}_{k}/\norm{\gls{kmd}_{k}}_{2}\) \tcp*{Mode normalization}
  \(\gls{kev}_{k} \gets (\ln \gls{kdev}_{k})/\Delta t\) \tcp*{Continuous-time DMD eigenvalues}
  \(\gls{kcoeff}\gets\argmin_{v} \sum_{n \subset [0,N-1]}\norm*{
\mat{X}[n] - \gls{kmd} \diag\begin{bmatrix} \gls{kdev}_{1}^{n}  & \gls{kdev}_{2}^{n}& \dots \end{bmatrix} v
  }^{2}_{2}\)
  \tcp*{\(\ell^{2}\) fit}
  \Output{\(\gls{kmd}_k, \gls{kdev}_{k},\gls{kev}_k, \gls{kcoeff}_{k}\)}
  \label{alg:modified-exact-DMD}
\end{algorithm}

\section{Analysis of the DMD spectrum and DMD modes}\label{sec:spectral-analysis}

The number of linearly independent modes that can be calculated from the full snapshot \(\mat{X}\) is its rank, which in our case is capped by the number of snapshots \(N=144\).
We now further want to extract a smaller number of modes in order to investigate them more deeply and to build a low-rank representation of data.

\subsection{SVD subspace selection}\label{sec:SVD_sub}

The first step in model reduction of DMD is the choice of parameter \(r\) which restricts the snapshot matrices \(\mat{X}_{1},\mat{X}_{2}\) to the \(r\)-dimensional subspace spanned by singular vectors of \(\mat{X}_{1}\).
Projections of \(X_{1}\) onto such subspaces amount to \(\ell^{2}\)-norm optimal approximations of \(X_{1}\) by \(r\)-ranked matrices (Eckart--Young theorem).
Since all subsequent calculations will be restricted to that subspace, the choice of \(r\) can dramatically affect the outcome of DMD\@.
There are no universal guidelines for choosing \(r\).
If the choice of \(r\) is analyzed separately from the DMD context in which we use it, then~\citep{Gavish2014} provides a threshold that removes iid, zero-mean additive noise,~\citep{Epps2019} reviews several different analytic and heuristic recipes for choosing the cutoff, while~\citep{Drmac2018} uses the numerical rank of the matrix (number of normalized singular values larger than a chosen threshold).
We remain agnostic with respect to the specific choice, and note that the dominant DMD modes are computed robustly as long as \(r\) is chosen away from the extreme values, for the following reasons.

Choosing \(r\) too small limits the number of linearly independent DMD modes that can be computed; typically, \(r\) should not be made smaller than the number of physically-relevant DMD modes, if there is prior information about them.
However, even if there is intuition about the physically-relevant modes, they may not be dominant in the \(\ell^{2}\) sense.
For this reason, one should choose \(r\) larger than expected, and perform the mode selection after all \(r\) DMD modes have been computed.

Choosing \(r\) too large can include into the subspace directions that are associated with noise (numerical or otherwise) into the computation.
DMD computations can be fragile with respect to the level of added noise, so it is prudent to remove at least some of the singular values in the tail.

\begin{figure}[htb]
  \centering
  \includegraphics[width=0.48\textwidth]{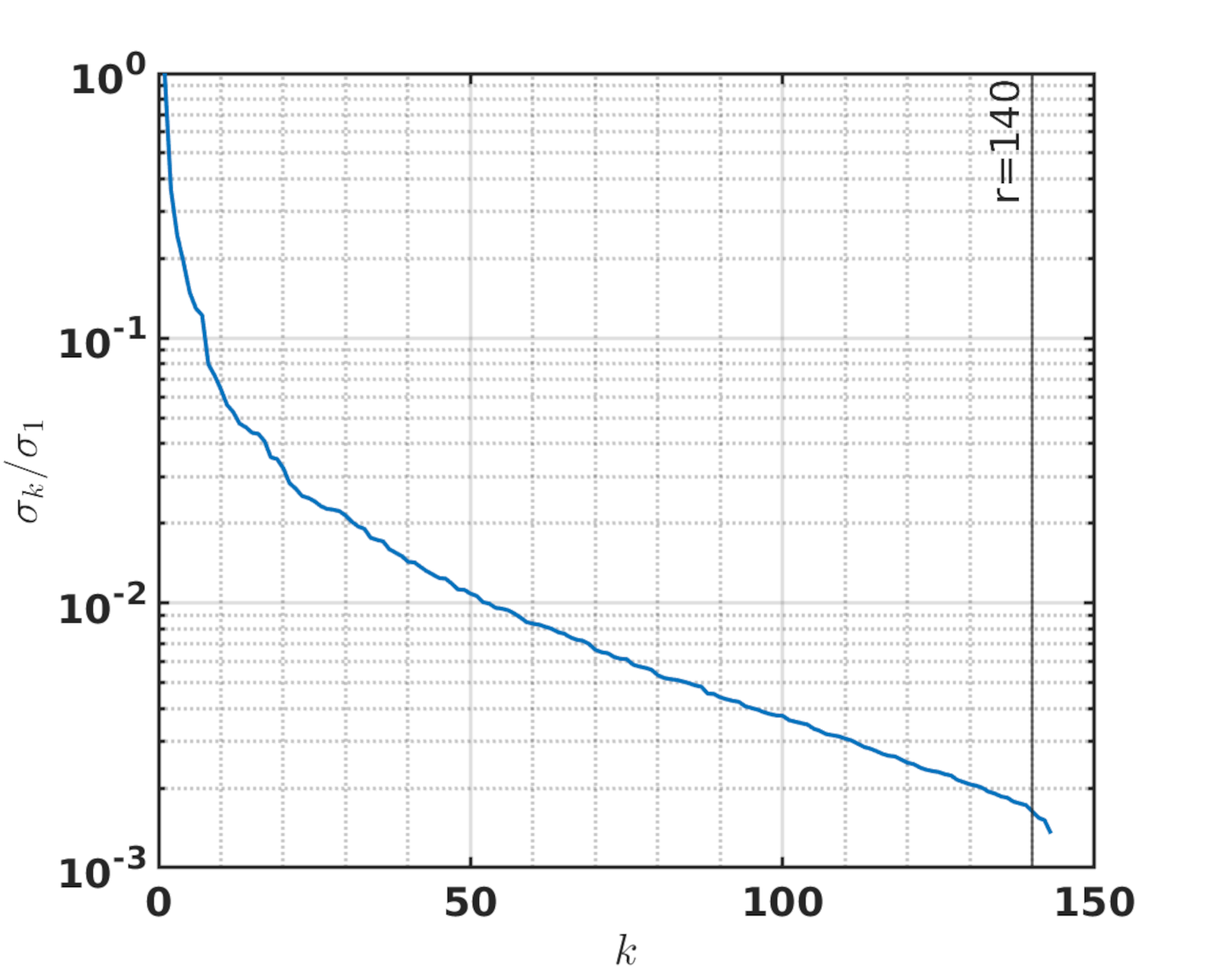}
  \caption[SVD values]{Singular values of the data matrix do not show significant separations of scale beyond the initial modes. The vertical line indicates the truncation value employed in our analysis.}
  \label{fig:svd-rank}
\end{figure}

For the remainder of the paper, we chose \(r=140\) as the cut off.
The singular values spectrum in Figure~\ref{fig:svd-rank} does not show particular spectral gaps, or leveling off, that may be associated with the presence of noise.
We have tested several additional choices of values \(r\) between 80 and 143 (the rank of the data matrix and number snapshots) without noticeable impact on the quality of the remainder of the analysis.

\subsection{Robustness of eigenvalues}\label{sec:robustness}

Algorithm~\ref{alg:modified-exact-DMD} serves to efficiently estimate eigenvalues and eigenvectors of the matrix \(\gls{dmdk}\) by effectively performing a linear regression on pairs of snapshots
\begin{equation}
    \gls{dmdk} := \argmin_{A} \norm{ \mat{X}_{2} - A\mat{X}_{1}}. \label{eq:DMD-matrix}
\end{equation}

If the number of snapshots is large, omitting a pair of snapshots would not distort resulting eigenvalues.
However, with a finite number of snapshots, the DMD eigenvalues do change slightly, but noticeably, even if a single snapshot is omitted.
When the expected frequencies of modes are well-separated, this likely does not pose a problem.
However, in the case of the analyzed data, periods of diurnal and semidiurnal tidal components used to force the model are tightly clustered around \(12\si{\hour}\) and \(24\si{\hour}\), differing sometimes only in third or fourth significant digit.
Consequently, the resulting small uncertainty in eigenvalues prevent us from clearly identifying a single DMD mode with a single tidal pattern.

To estimate the robustness of the linear regression with respect to the amount of data given, we perform the leave-one-out analysis.
After computing the DMD eigenvalues and modes for the full data set, we repeat the computation of eigenvalues 30 times while erasing at random \(i\)th columns \(\mat{X}_{1}[:,i], \mat{X}_{2}[:,i]\) from the pair of matrices.

The resulting collection of eigenvalue spectra is used in a Kernel Density Estimates (KDE) of the distribution of eigenvalues.
A set of \(30 \times r\) eigenvalues \(\hat{\gls{kdev}}_{i}\), where \(r=140\) is the dimension of the SVD subspace, is converted to a continuous density on the complex plane \(z \in \mathbb{C}\) via
\begin{equation}
  \label{eq:KDE-eigenvalues}
  d_{h}(z) \propto \sum_{k} \exp \frac{-\abs{z - \hat{\gls{kdev}}_{k}}^{2}}{h^{2}},
\end{equation}
where \(h\) is the bandwidth parameter that controls the width of each Gaussian superimposed at locations of eigenvalues, and the density is normalized to unit mass.
Choosing \(h\) wider or narrower produces a smoother or more localized density estimate.
Eigenvalues that do not move significantly between leave-one-out experiments result produce density that is concentrated over small subsets of complex plane, while high-variability eigenvalues result in a smeared-out parts of the density.

\begin{figure*}[htb]
  \centering
  \includegraphics[width=\textwidth]{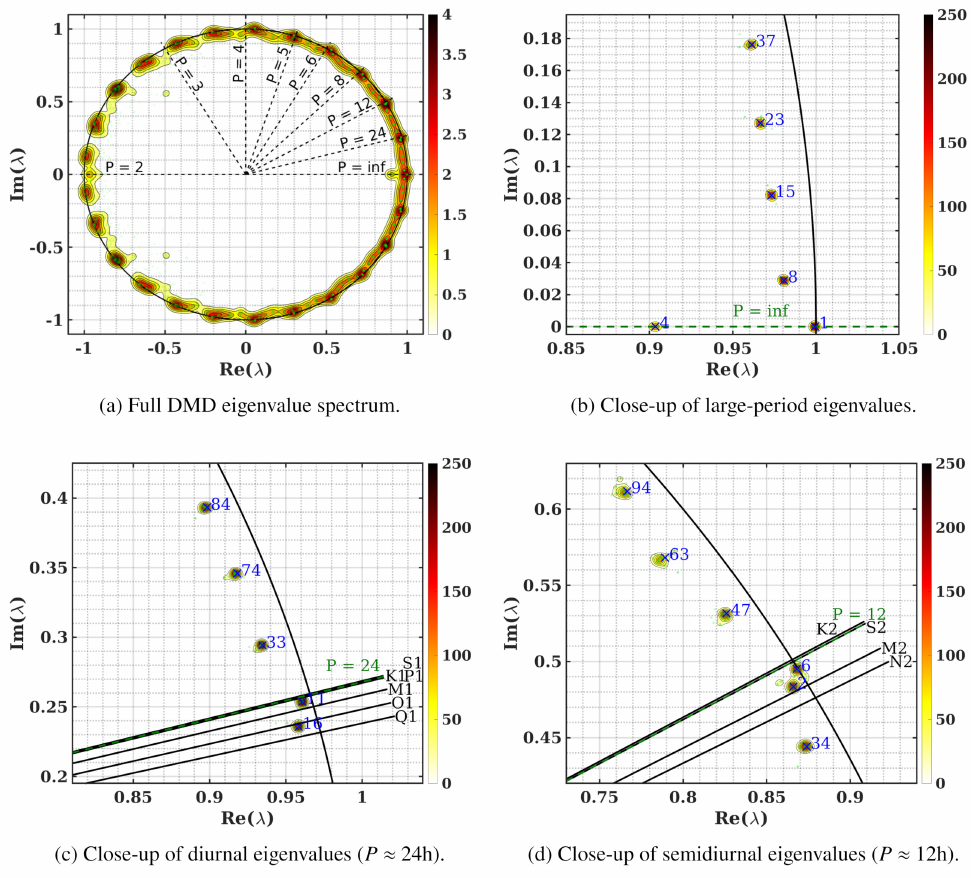}
\caption[Leave-one-out analysis]{The spectrum of DMD discrete-time eigenvalues~\(\gls{kdev}_k\). Unit circle separates decaying modes (inside) from growing modes (outside). Diagonal lines indicate frequencies corresponding to tidal components used in forcing.
Colored contours show Gaussian kernel density estimates (KDE) based on leave-one-out analysis. Kernel bandwidth is \(h=2.5\times 10^{-2}\) in panel (a) and \(h=2.0\times 10^{-3}\) in zoomed-in panels (b)-(d). Blue crosses indicate eigenvalues computed using all the snapshots with numbers corresponding to the RMS rank of each mode as reported in Table~\ref{tab:eigenvalues}}.\label{fig:leave-one-out}
\end{figure*}

In Figure~\ref{fig:leave-one-out}, eigenvalues aggregate in well-defined clusters around frequencies associated with tidal components used to force the model. 
 An exception is the robust cluster around \(P \approx 5 \si{\hour}\) which is not a harmonic of any input frequency.
 As expected, shorter time periods, where only a few snapshots are available per period, are less robustly computed.
 Eigenvalues close to the unit circle, e.g., as seen in panel Figure~\ref{fig:leave-one-out}(d) will commonly result in leave-one-out perturbations that straddle the unit circle, that is result in slightly growing or slightly decaying modes; this analysis shows that such behavior can be the result of finite-length data set, rather than physical growth or decay.

In addition to visualization, the density \(d_{h}(z)\) estimated from leave-one-out eigenvalues \(\hat{\gls{kdev}}_{k}\)can be evaluated at locations of eigenvalues \(\gls{kdev}_{n}\) computed for the full set of snapshots.
The values \(d_{h}(\gls{kdev}_{n})\) provide an estimate of robustness of eigenvalues \(\gls{kdev}_{n}\)to perturbations in available data; they are reported in the summary Table~\ref{tab:eigenvalues}.

\begin{table*}[hbt]
  \centering
  \begin{minipage}{.495\textwidth}
    \centering
    \resizebox{\textwidth}{!}{
      \begin{tabular}{@{\extracolsep{\fill}}ccccccc}
        \topline 
        \thead{Index} & \thead{Cluster} & \thead{Period} & \thead{Halving (-) /\\ Doubling (+) Time} & \thead{RMS \\ \(\gls{meanL2}_{n}\)} & \thead{\(\gls{meanL2W}_{n}\)} & \thead{Robustness\\\(d_{h}(\gls{kev}_{n})\)} \\
        \midline
   1 &        1 &    inf & -1078.97 & 157.02 &   3.24 &   282.78 \\
   2 &        3 &  12.34 &   -84.48 &  75.46 &   4.00 &   181.47 \\
   4 &        1 &    inf &    -6.84 &  74.58 &   1.72 &   163.53 \\
   6 &        3 &  12.12 & -1029.16 &  69.55 &   3.79 &   190.64 \\
   8 &        1 & 212.95 &   -36.32 &  56.86 &   2.84 &   257.16 \\
  10 &       NaN &   2.80 &    -3.26 &  36.87 &   0.15 &     0.00 \\
  11 &        2 &  24.34 &  -112.00 &  33.99 &   3.50 &   236.18 \\
  13 &       13 &   2.00 &    -3.59 &  27.67 &   0.08 &     0.00 \\
  15 &        1 &  74.50 &   -29.33 &  25.15 &   2.34 &   230.53 \\
  16 &        2 &  26.04 &   -53.05 &  19.28 &   3.40 &   217.63 \\
  19 &       10 &   2.81 &    -7.39 &  16.41 &   0.45 &     0.00 \\
  21 &       13 &   2.14 &    -4.58 &  16.29 &   0.20 &     0.01 \\
  23 &        1 &  48.02 &   -27.34 &  15.04 &   2.41 &   197.61 \\
  24 &        5 &   6.19 &    90.93 &  11.11 &   3.08 &   199.52 \\
  26 &       11 &   2.57 &    -7.01 &  10.77 &   0.61 &     0.00 \\
  28 &        9 &   2.97 &    -8.75 &   9.75 &   0.61 &     2.92 \\
  30 &        4 &   8.26 & -2531.05 &   9.29 &   3.29 &   210.09 \\
  33 &        2 &  20.60 &   -34.22 &   8.85 &   2.35 &   212.79 \\
  34 &        3 &  13.36 &   -34.92 &   8.84 &   2.65 &   186.66 \\
  37 &        1 &  34.67 &   -30.36 &   8.51 &   2.19 &   188.51 \\
  39 &        6 &   4.95 &   -36.89 &   8.20 &   3.31 &    37.57 \\
  41 &       12 &   2.29 &    -8.28 &   7.60 &   0.52 &     0.00 \\
  43 &        6 &   4.95 &  -943.33 &   7.44 &   3.30 &   130.84 \\
  45 &        4 &   7.92 &   -47.62 &   6.99 &   2.75 &   102.78 \\
  47 &        3 &  10.99 &   -38.63 &   6.56 &   2.25 &   157.12 \\
  49 &        4 &   7.43 &   -22.31 &   6.14 &   2.08 &    33.50 \\
  50 &       12 &   2.38 &   -10.92 &   6.05 &   0.70 &     4.73 \\
  53 &        5 &   5.97 &   -73.14 &   5.83 &   2.46 &    65.46 \\
  55 &        5 &   6.99 &   -18.17 &   5.41 &   2.20 &    19.44 \\
  56 &        7 &   4.15 &   140.13 &   5.11 &   3.16 &   107.07 \\
  58 &        9 &   3.32 &   -10.05 &   5.08 &   0.94 &     2.27 \\
  60 &        5 &   6.34 &   -52.22 &   5.02 &   2.89 &   107.70 \\
  63 &        3 &  10.07 &   -24.71 &   4.97 &   2.22 &    62.55 \\
  65 &        6 &   4.75 &   -35.73 &   4.90 &   2.17 &    12.03 \\
  66 &       13 &   2.00 &   -13.31 &   4.88 &   1.61 &     8.35 \\
  67 &        4 &   8.56 &   -44.77 &   4.83 &   2.88 &   141.89 \\
  70 &       10 &  2.72 &   -25.92 &  4.65 &   1.77 &     5.55 \\
        \botline
      \end{tabular}
    }
  \end{minipage}
  \begin{minipage}{.495\textwidth}
    \centering
    \resizebox{\textwidth}{!}{
      \begin{tabular}{@{\extracolsep{\fill}}ccccccc}
        \topline 
        \thead{Index} & \thead{Cluster} & \thead{Period} & \thead{Halving (-) /\\ Doubling (+) Time} & \thead{RMS \\ \(\gls{meanL2}_{n}\)} & \thead{\(\gls{meanL2W}_{n}\)} & \thead{Robustness\\\(d_{h}(\gls{kev}_{n})\)} \\
        \midline
  71 &        7 &  3.94 &   -15.78 &  4.49 &   1.67 &     3.78 \\
  74 &        2 & 17.44 &   -36.43 &  4.43 &   2.20 &   179.05 \\
  76 &        7 &  4.06 &  -141.09 &  4.11 &   2.47 &    35.12 \\
  78 &        5 &  5.70 &   -27.64 &  3.79 &   2.09 &    15.80 \\
  80 &        7 &  4.54 &   -18.55 &  3.78 &   1.75 &     6.04 \\
  82 &        8 &  3.47 &   -38.53 &  3.66 &   2.17 &    10.31 \\
  84 &        2 & 15.22 &   -35.50 &  3.49 &   2.17 &   184.71 \\
  85 &        9 &  3.07 &  -107.53 &  3.30 &   2.53 &    10.21 \\
  88 &        8 &  3.53 &   212.17 &  3.05 &   2.76 &   122.36 \\
  90 &        5 &  6.69 &   -24.35 &  2.81 &   2.32 &    41.60 \\
  92 &        9 &  3.11 &  -249.05 &  2.70 &   3.33 &    32.26 \\
  94 &        3 &  9.33 &   -35.46 &  2.44 &   2.13 &    99.96 \\
  95 &       12 &  2.23 &   -49.59 &  2.38 &   3.00 &    13.64 \\
  98 &       10 &  2.75 &  6512.72 &  2.27 &   2.79 &    49.60 \\
 100 &       11 &  2.47 &  -557.29 &  1.99 &   2.55 &    18.20 \\
 101 &        8 &  3.68 &   -21.23 &  1.87 &   2.15 &     0.31 \\
 104 &        8 &  3.60 &   -47.45 &  1.57 &   2.92 &    22.00 \\
 106 &        7 &  4.38 &   -31.32 &  1.55 &   1.92 &     4.40 \\
 107 &        6 &  5.13 &   -48.67 &  1.52 &   2.44 &    65.01 \\
 109 &        8 &  3.75 &   -35.62 &  1.42 &   1.91 &     3.94 \\
 111 &        7 &  4.23 &   -50.03 &  1.38 &   2.74 &    37.96 \\
 114 &        6 &  5.41 &   -41.12 &  1.36 &   1.83 &    19.34 \\
 116 &       11 &  2.51 &   -58.72 &  1.32 &   4.19 &     0.02 \\
 117 &       12 &  2.26 & -1149.21 &  1.13 &   3.07 &    38.18 \\
 120 &        9 &  3.16 &  -103.37 &  1.07 &   2.74 &    55.52 \\
 121 &        9 &  3.21 &   -66.38 &  1.06 &   2.06 &     2.77 \\
 123 &       13 &  2.06 &  -614.10 &  0.94 &   2.05 &    31.39 \\
 125 &       11 &  2.42 &   -20.55 &  0.89 &   0.66 &     0.00 \\
 128 &       13 &  2.11 &  -386.21 &  0.83 &   1.95 &     1.20 \\
 129 &       11 &  2.53 &  -102.69 &  0.79 &   3.49 &    27.99 \\
 132 &       12 &  2.30 & -2059.91 &  0.75 &   1.98 &    60.64 \\
 133 &       13 &  2.09 &  -251.43 &  0.72 &   2.84 &     1.12 \\
 136 &       13 &  2.14 &   -21.67 &  0.67 &   0.81 &     2.09 \\
 138 &       10 &  2.78 &   -28.79 &  0.64 &   2.64 &     0.00 \\
 139 &       10 &  2.81 &  -124.00 &  0.51 &   2.90 &    38.63 \\
 & & & & & & \\
        \botline
      \end{tabular}
    }
  \end{minipage}
  \caption{Properties of DMD modes.
Only one conjugate is listed for each complex mode pair.
Robustness is estimated using kernel width \(h=2\times 10^{-3}\).
Clusters were assigned based on level sets of kernel density estimates with \(h=2.5\times 10^{-2}\); label ``NaN'' indicates modes that  do not belong to a well-defined cluster.\label{tab:eigenvalues}}
\end{table*}

\subsection{Contributions of individual modes}\label{sec:mode-ranking}

There is no \emph{a priori} ranking of DMD modes that can be used to easily define a criterion for the ranking of modes by their ``importance''.
The \(\ell^{2}\) norm of each mode, \(\abs{\gls{kcoeff}_{k}}\), is one natural model-agnostic candidate, but decay rates (to identify persistent modes), \(\ell^{2}\) contributions to a subset of physical variables, or other characterizations of modes (see~\citep{bollt2021geometric}) may be equally or even more important depending on the model analyzed using DMD\@.
We here showcase several such possibilities in the context of ocean dynamics.

Magnitudes of coefficients \(\abs{\gls{kcoeff}_{k}}\) used to assemble the model of the evolution reflect how much each mode contributes to the initial condition~\eqref{eq:initial-condition}.
Unfortunately, such comparison will favor modes that are decaying, as their magnitude is initially large and then decreases in time.

To account for this, we calculate the Root Mean Square (RMS) norm of each mode over time
  \begin{align}
\label{eq:mean-l2-contribution}
\gls{meanL2}_{i}^{2} &=\frac{1}{T} \int_{0}^{T} \norm{\gls{kcoeff}_i e^{\gls{kev}_i t} \gls{kmd}_{i} }^2 dt = \abs{\gls{kcoeff}_{i}}^{2} \frac{e^{2 \gls{kdr}_i T}-1}{2 \gls{kdr}_i T} \\
\gls{meanL2}_{i} &= \abs{\gls{kcoeff}_{i}} (1 + \gls{kdr}_{i} T/2 + \dots),
\end{align}
where \(\gls{kdr}_{i} = \Re \gls{kev}_{i}\).
The expansion demonstrates that the correction will reduce the significance of decaying (\(\gls{kdr}_{i} < 0\)) modes and boost the significance of growing modes, but only when \(\abs{\gls{kdr}_{i}T}\) is close to \(1\) or larger.

\begin{figure*}[htb]
  \centering
  \includegraphics[width=\textwidth]{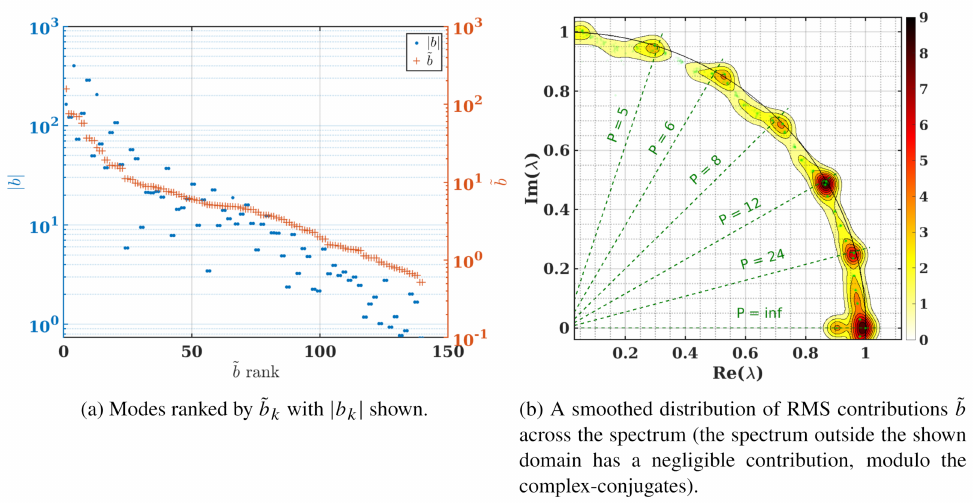}
\caption{The distribution of RMS \(\gls{meanL2}_{k}\) across the eigenvalue spectrum and comparison with \(\abs{\gls{kcoeff}_{k}}\).}\label{fig:RMS-L2norm}
\end{figure*}

Figure~\ref{fig:RMS-L2norm} shows the difference in ranking the modes using mode contributions to the initial condition \(\gls{kcoeff}\) vs.~RMS contributions of modes.
Majority of modes would be similarly ranked toward the beginning or end of the spread regardless of the choice, but for some the difference is significant.
For purposes of visualization, the concentration of RMS contributions across the spectrum of DMD eigenvalues is represented by the heat map of kernel density estimate of eigenvalues in Figure~\ref{fig:RMS-L2norm}\textbf{(b)}.
Each eigenvalue has been weighted by its RMS coefficient \(\gls{meanL2}\) to yield the function
\begin{equation}
  \label{eq:energy-density}
  \gls{meanL2}(z) \coloneqq  \sum_{k=1}^{r} \gls{meanL2}_{k} \exp \frac{-\abs{z - \gls{kdev}_{k}}^{2}}{h^{2}},
  \end{equation}
  where kernel bandwith \(h=2 \times 10^{-3}\) was chosen.
  It is clear that semidiurnal and diurnal eigenvalues dominate in terms of RMS contribution, followed by long-period and terdiurnal (\(P \approx 8 \si{hour}\)) eigenvalues.
  Notice also the \(P \approx 5 \si{hour}\) group of modes which does not correspond to either principal tides or their harmonics (``overtides'').

\begin{figure*}[htb]
  \centering
  \includegraphics[width=\textwidth]{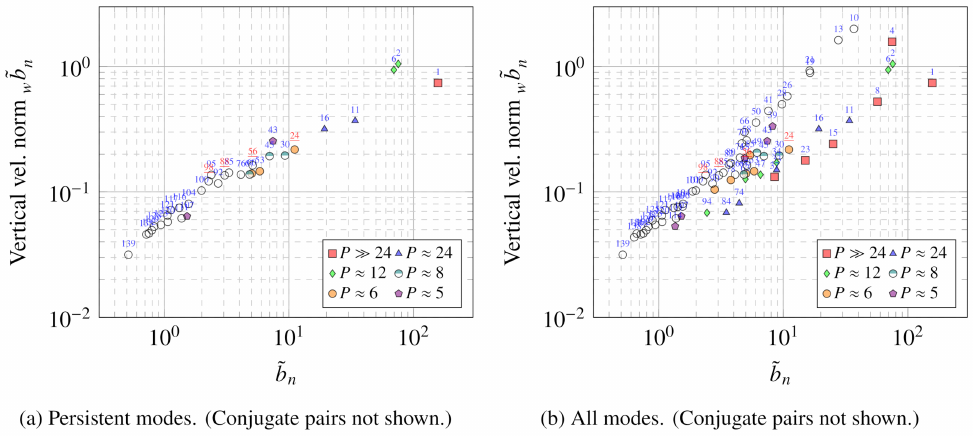}
\caption{Comparison between RMS of the total and the vertical velocity norms of modes.
  Indexes of modes correspond to Table~\ref{tab:eigenvalues}; marker shapes are used to indicate memberships in primary clusters of interest classified according to their oscillation periods \(P\); underlined labels correspond to growing modes.
Left panel omits \emph{non-persistent} modes ~\eqref{eq:ten-percent-decay}.}
\label{fig:total_vs_vertical}
\end{figure*}

\begin{figure*}[htb]
  \centering
  \includegraphics[width=\textwidth]{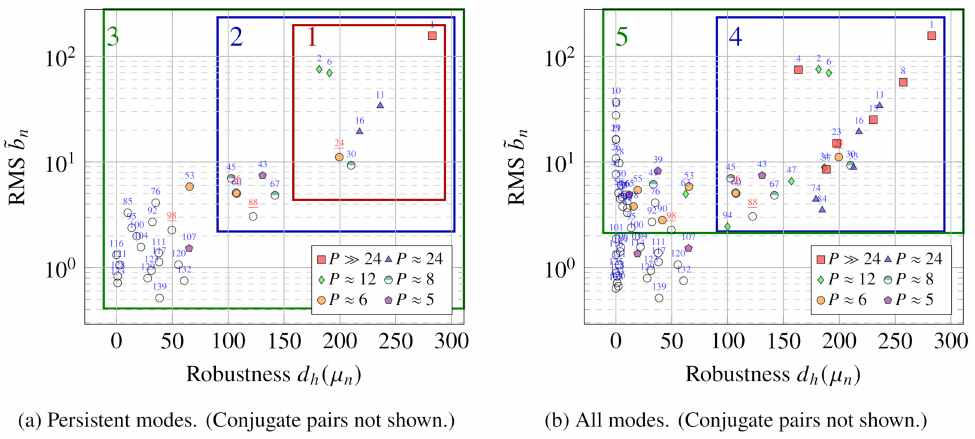}
  \caption{Arrangement of DMD modes with respect to RMS \(\gls{meanL2}_{k}\) and robustness \(d_{h}(\mu_{k})\)  with kernel width \(h=2\times10^{-3}\).
  Indexes of modes correspond to Table~\ref{tab:eigenvalues}; marker shapes are used to indicate memberships in primary clusters of interest classified according to their oscillation periods \(P\); underlined labels correspond to growing modes.
  Left panel omits \emph{non-persistent} modes~\eqref{eq:ten-percent-decay}.
  Rectangles indicates modes collected into five different reduced order models (see Section~\ref{sec:low-rank}), comprising, respectively, 13, 25, 65, 46 and 98 modes, amounting to approximately 9, 17, 45, 32, and 69 percent of the rank of the data matrix.}   \label{fig:robustness-energy-trade-off}
\end{figure*}

Model-agnostic methods of ranking the modes, such  as norm-contributions \(\gls{kcoeff}_{k}\), \(\gls{meanL2}_{k}\) would in practice be supplemented by model-specific considerations.
For example, in the particular case of tides one could keep harmonics of fundamental tides (``overtides'') as the part of the model, even when their RMS is not very high.
Similar discussions along these lines for EOF/POD modeling can be found in~\citep{tissot2014model,Monahan2009,Crommelin2004}.
An example of discrimination based on the decay rates is shown by the distinction between left and right panels in Figure~\ref{fig:total_vs_vertical} and Figure~\ref{fig:robustness-energy-trade-off} where modes are deemed non-persistent if their initial magnitudes decays to less than \(10\%\) of the initial magnitude, i.e., if
\begin{align}
  \label{eq:ten-percent-decay}
  e^{\gls{kdr} T } < 0.1, &\text{ for } T = 143 \si{hour}
  \shortintertext{which translates into half-life times}
  -43 < &T_{1/2} < 0.
\end{align}
One purpose of such selection of variables would be to discard the initial transient, for example, in simulations that are not initialized on the attractor of dynamics.

The vertical velocity in the simulated ocean is by, roughly, two orders of magnitude smaller than the total velocity; therefore comparing modes by total velocity potentially obscures how the vertical transport is distributed across the modes.
To evaluate what modes are important for the vertical transport we additionally calculate the \(\ell^{2}\) norm of mode elements corresponding to vertical velocities \(\gls{L2W}_{k}\), and the associated RMS contribution \(\gls{meanL2W}_{k}\), analogously to~\eqref{eq:mean-l2-contribution}.
Figure~\ref{fig:total_vs_vertical} shows that the magnitudes of total and vertical velocities correlate.
In both cases, the correlation roughly linear, indicating that the modes with strongest vertical activity in most cases rank highly in terms of the norm of the full velocity field.
Further exploration of the correlations is beyond the scope of this paper.

Finally,  Figure~\ref{fig:robustness-energy-trade-off} demonstrates correlations between robustness of eigenvalues, as estimated by \(d_{h}(\gls{kev}_{k})\) estimate~\eqref{eq:KDE-eigenvalues}, and the contributions, as measured by \(\gls{kcoeff}_{k}\) or \(\gls{meanL2}_{k}\).
  The robust, high-contributing modes would be located in the upper right corner of each of the graphs, and are the top candidates for building any low-rank model of the data.
  Modes that score well on one, but not both axes are in top left, and bottom right, and depending on the purpose of a low-rank model, one or both of these groups could be included.
  Modes close to the origin of the coordinate system are neither robustly estimated nor norm-significant and would be likely discarded.

  The structure of spatial profiles \(\Phi\) can be investigated in other ways; for example,~\citep{bollt2021geometric} clarifies the importance of the total variation norm (the norm of the spatial gradient) of \(\Phi\) in identification of primary modes.
  We leave the pursuit of such analyses to a follow-up publication.

\subsection{Low-rank model of simulated data}
\label{sec:low-rank}

To compare several different reduced-order models (ROMs) of data we retained the modes as indicated by rectangles in Figure~\ref{fig:robustness-energy-trade-off}.
The rectangles were chosen based on natural separation of modes in groups on these graphs, although certainly more sophisticated criteria could be formulated, including various clustering techniques and other data mining methods.

\begin{figure*}[htb]
  \centering
  \includegraphics[width=\textwidth]{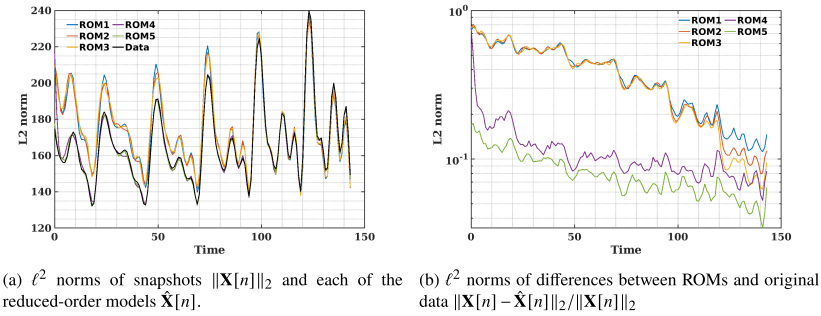}
  \caption{Quality of approximation of the original data set \(\mat{X}[n]\) using each of the reduced-order modes \(\hat{\mat{X}}[n]\) corresponding to selections indicated in Figure~\ref{fig:robustness-energy-trade-off}.}
\label{fig:ROM-l2-norms}
\end{figure*}

The dimensions of ROMs are respectively 13, 25, 65, 46 and 98 modes, amounting to approximately 9, 17, 45, 32, and 69 percent of the rank of the data.
ROMs 1--3 include only persistent modes, in the sense of~\eqref{eq:ten-percent-decay}, while ROMs 4 and 5 include all modes within the indicated rectangles.

Figure~\ref{fig:ROM-l2-norms} shows norms of low-rank approximations (reduced-order models, ROMs) to the data and norms of differences between data \(\mat{X}\) and each ROM \(\hat{\mat{X}}\).
It is clear that all ROMs, despite the wide range of their orders, capably reproduce bulk properties (norm) of the data in the latter part of the time domain.
The clear difference between ROMs 1--3 and 4-5 is in the transient part; this is expected, since ROMs 1--3 intentionally used only the modes that have a relatively-slow decay rate, or are growing.
Within the two classes of ROMs (persistent-only and full), the differences are minor even though the dimension of ROMs vary considerably within each class.

\section{Connections between modes and oceanographic features}
\label{sec:connections-oceans}

In the following subsections, we highlight connections between oceanographic features and modes, or groups of modes, isolated by the DMD analysis.
In particular, we focus on the most persistent dominant modes, enclosed in ``ROM 1'' box in Figure~\ref{fig:robustness-energy-trade-off}, and listed in Table~\ref{tab:ROM1-modes}, as these are the most reliably computed features that are observable over the duration of the simulation.

\begin{table*}[t]
\begin{center}
\begin{tabular}{cccccccc}
\topline
Index & Cluster & Period & Halving ($-$)/Doubling ($+$) Time & RMS \(\gls{meanL2}_{n}\) & \(\gls{meanL2W}_{n}\) & Robustness $d_h(\lambda_n)$ \\
\midline
1 & 1 & Inf & $-1078.97$ & $157.02$ & $0.74$ & $282.78$ \\
2 & 3 & $12.34$ & $-84.48$ & $75.46$ & $1.05$ & $181.47$ \\
6 & 3 & $12.12$ & $-1029.16$ & $69.55$ & $0.94$ & $190.64$ \\
11 & 2 & $24.34$ & $-112.00$ & $33.99$ & $0.37$ & $236.18$ \\
16 & 2 & $26.04$ & $-53.05$ & $19.28$ & $0.32$ & $217.63$ \\
24 & 5 & $6.19$ & $+90.93$ & $11.11$ & $0.22$ & $199.52$ \\
30 & 4 & $8.26$ & $-2531.05$ & $9.29$ & $0.20$ & $210.09$ \\
\botline
\end{tabular}
\end{center}
\caption{Properties of DMD modes included in the ``ROM 1'' reduced order model. 
Only one conjugate is listed for each complex mode pair.
Robustness is estimated using kernel width \(h=2\times 10^{-3}\).
Clusters were assigned based on level sets of kernel density estimates with \(h=2.5\times 10^{-2}\);
Listed modes are a subset of the modes in Table~\ref{tab:eigenvalues}.} \label{tab:ROM1-modes}
\end{table*}

\subsection{Western Alboran Gyre and secondary gyres}\label{sec:WA_SG}
\begin{figure*}[htb]
  \centering
  \includegraphics[width=\textwidth]{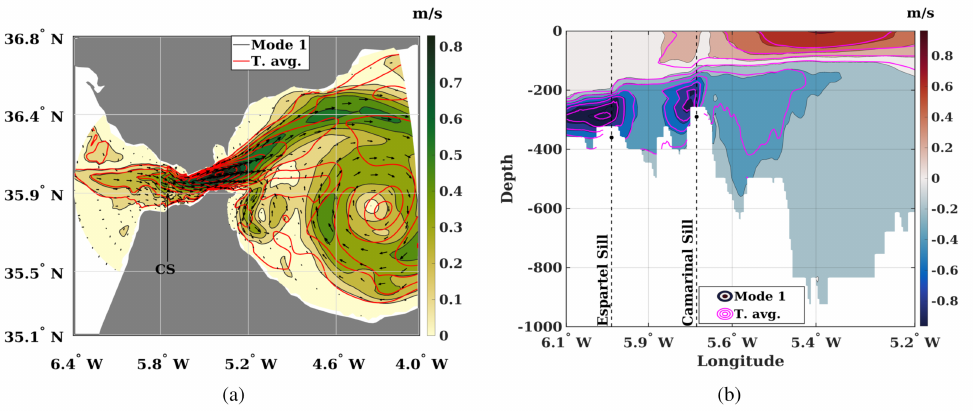}
  \caption{(a) Horizontal velocity field; heatmap is the magnitude of the horizontal velocity vector. The contours in red show the data, $S_{\overline{U}}$, before applying DMD. (b) Mode 1 zonal velocity ($U_x$). Contours (in magenta) are the time-average of the same quantity.}
\label{fig:11_R5}
\end{figure*}

Mode 1 is an aperiodic, near-constant mode that captures the largest proportion of the time-averaged norm of data. For this reason, we expect it to show flow features similar to the time-averaged data, displayed in Figure~\ref{fig:1_R1}. Figure~\ref{fig:11_R5} show that this is largely the case: the filled-in contours of Mode 1 align with overlaid contours of the time-averaged data. The agreement pertains to both the horizontal circulation, the vertical structure of the exchange flow within the strait, and even the qualitative structure of stationary internal waves created by flow over topographic features.

The only major difference between time-averaged data and the Mode 1 is in the region immediately southwest of the Strait, off Ceuta. In this region, the outer edge of the WAG reverses direction and continues to flow down African coast. Mode 1 shows that a secondary cyclonic gyre sets up as a result of the flow reversal in this region. The time-averaged velocity field does show increased flow intensity in that region, but it is not apparent from the flow direction in Figure~\ref{fig:11_R5}a that the flow is organized in a gyre, rather than a hairpin pattern. 

\subsection{Sill flow and internal waves }\label{sec:sl_FW}

\begin{figure*}[htb]
  \centering
  \includegraphics[width=\textwidth]{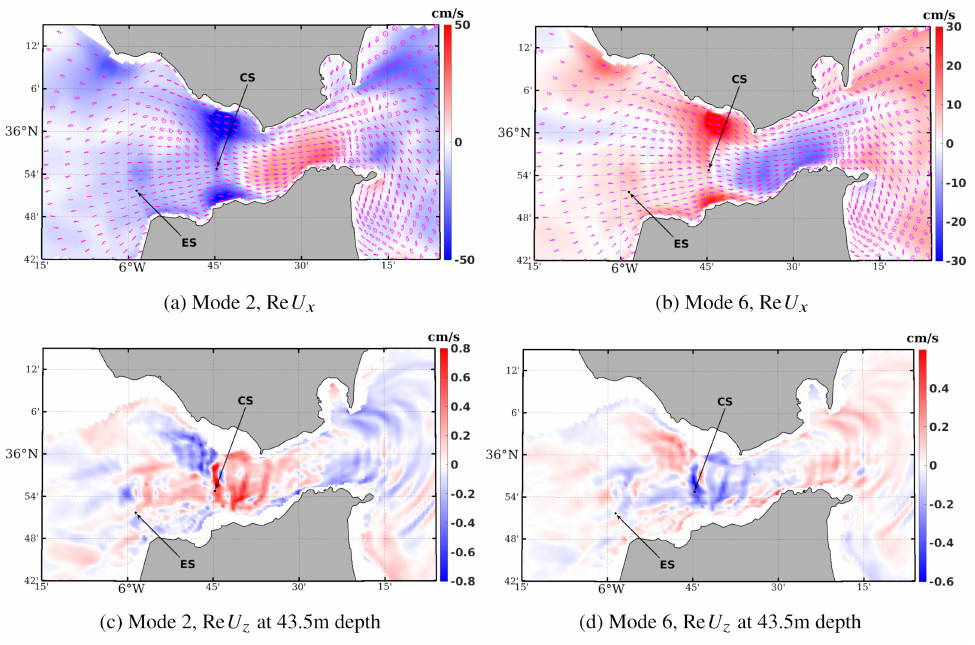}
  \caption{(a)-(b) Surface tidal ellipses and zonal velocity ($U_x$) of modes 2 and 6. (c)-(d) Subsurface vertical velocity ($U_z$) of the same modes, 2 and 6, at 43.5 m depth.}
  \label{fig:R6}
\end{figure*}

Modes 2 and 6 in Table~\ref{tab:ROM1-modes}, with respective periods of 12.34h and 12.12h are both approximately semidiurnal and therefore match the dominant period of variability in the strait. These modes are associated with a modulation of the horizontal velocity that is predominantly longitudinal, as evidenced by the tidal ellipses (Figure~\ref{fig:R6}a,b). The surface velocity field for modes 2 and 6 are horizontally divergent and convergent, respectively, immediately east of the CS, and the associated subsurface vertical velocities show dynamically consistent upwelling and downwelling patterns in this region (Figure~\ref{fig:R6}c,d). Therefore, modes 2 and 6 capture the large vertical motions east of the sill that can ultimately lead to the formation of the approach control when the interface (S=37.5) here becomes sufficiently shallow (Figure~\ref{fig:2_R2}c,d). These modes also display short-wave patterns at the eastern exit of the strait (Figure~\ref{fig:R6}c and Figure~\ref{fig:R6}d), likely associated with the internal bore. 

\begin{figure*}[htb]
  \centering
  \includegraphics[width=\textwidth]{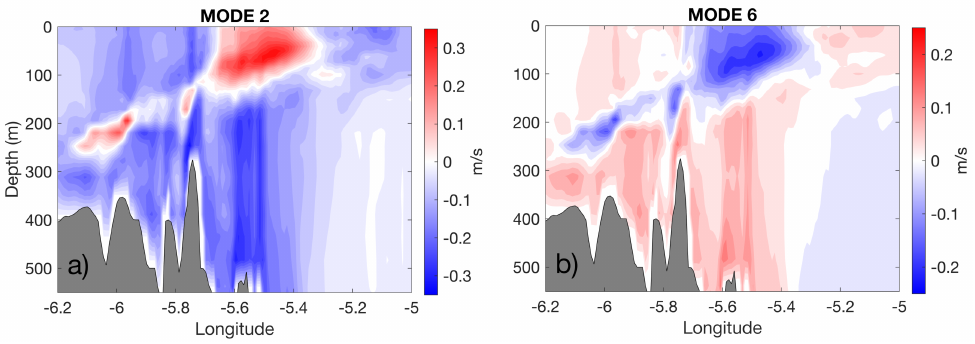}
  \caption{Along-strait cross section of zonal velocity $U_x$ for mode 2 (a) and mode 6 (b).}
  \label{fig:R7}
\end{figure*}

\begin{figure*}[htb]
  \centering
  \includegraphics[width=\textwidth]{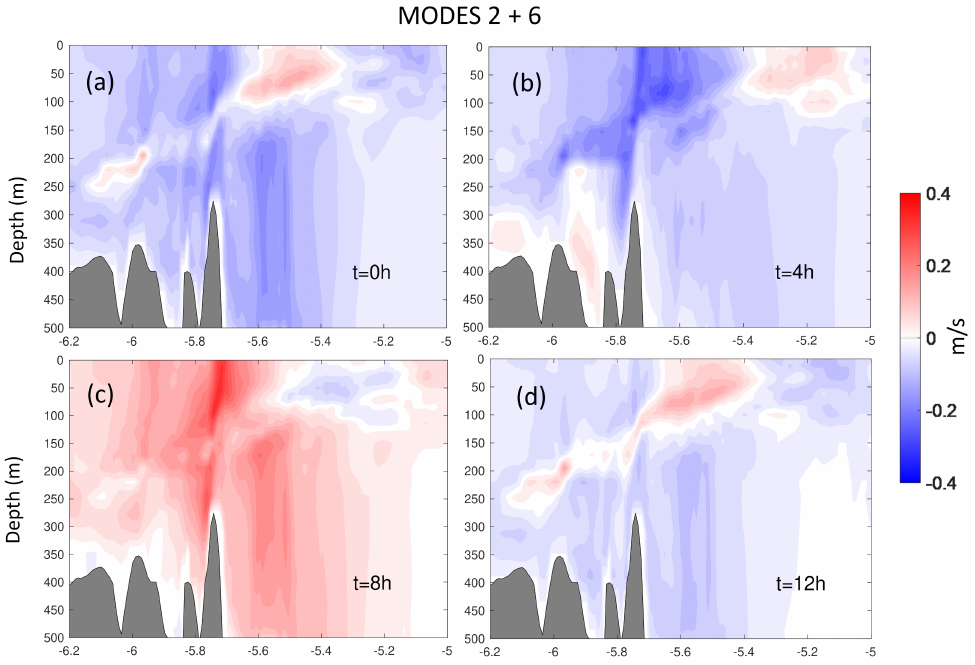}
  \caption{((a)-(d) Zonal velocity of modes 2 and 6 (superposition of the two) over a semidiurnal tidal cycle. Current-countercurrent oscillations east of CS are due to baroclinic tides propagating towards the Alboran Sea.}
  \label{fig:R8}
\end{figure*}

\begin{figure*}[htb]
  \centering
  \includegraphics[width=\textwidth]{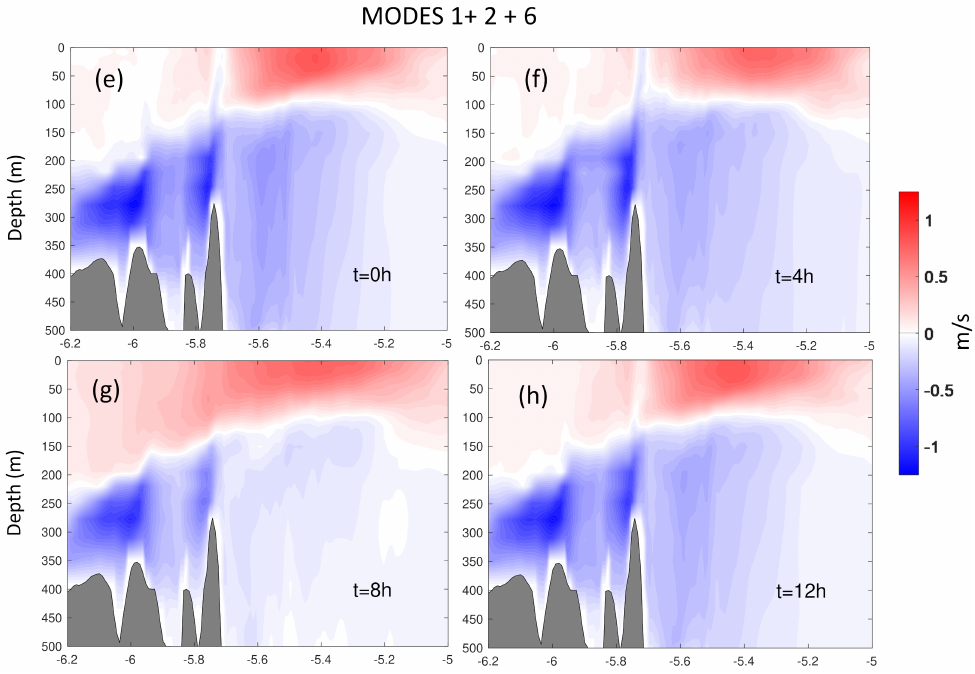}
  \caption{(e)-(f) Zonal velocity of modes 1, 2, and 6 (superposition of the three) over a semidiurnal tidal cycle.}
  \label{fig:R9}
\end{figure*}

Further inspection of Modes 2 and 6 reveals that horizontal currents have distinct vertical structure at either side of the CS, especially for Mode 2 (Figure~\ref{fig:R7}a); west of CS and over CS itself tides are mostly barotropic (Tsimplis, 2000), whereas east of CS tides have a stronger baroclinic structure, in agreement with observations (e.g., \cite{lafuente2000tide}). The horizontal velocity fields of Mode 2 and 6 exhibit a first-mode vertical baroclinic structure (i.e. a single zero crossing of horizontal velocity) east of CS, with both modes also displaying isolated signs of second-mode vertical structure (two zero crossings) west of the sill, especially Mode 6 (Figure~\ref{fig:R7}b). The horizontal velocity field of combined modes 2 and 6 over a tidal cycle shows that these modes account for the eastward radiation of first-mode internal tides from the CS and the strong barotropic tidal flow that periodically blocks, and even reverses, the Atlantic and Mediterranean currents over it (Figure~\ref{fig:R8}b,c). From the superposition of mode 1 and the semidiurnal modes 2 and 6, it can be shown that the Atlantic current reverses over the CS at t=4 h (Figure~\ref{fig:R9}f),while the Mediterranean current becomes nearly blocked at t=8h (Figure~\ref{fig:R9}g). These two instants correspond to flood and ebb tidal phases, respectively (i.e., westward and eastward barotropic tidal flow). 

\subsection{Meanders in the Atlantic jet}\label{sec:M_AJ}

\begin{figure*}[htb]
  \centering
  \includegraphics[width=\textwidth]{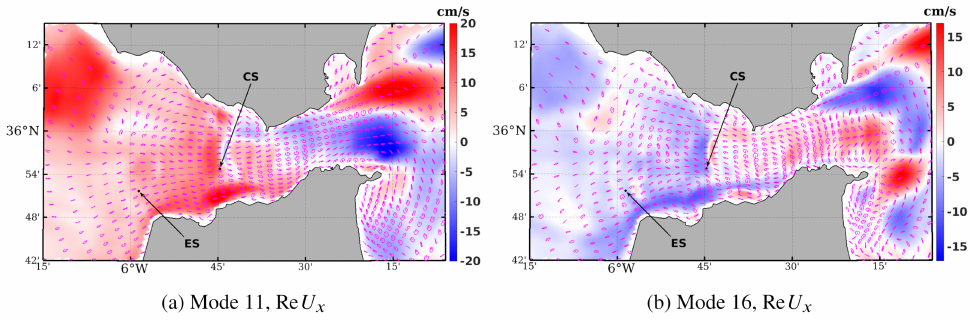}
\caption{Zonal surface velocity (Ux) of the diurnal modes 11 and 16.}
\label{fig:R9_1}
\end{figure*}

\begin{figure*}[htb]
  \centering
  \includegraphics[width=\textwidth]{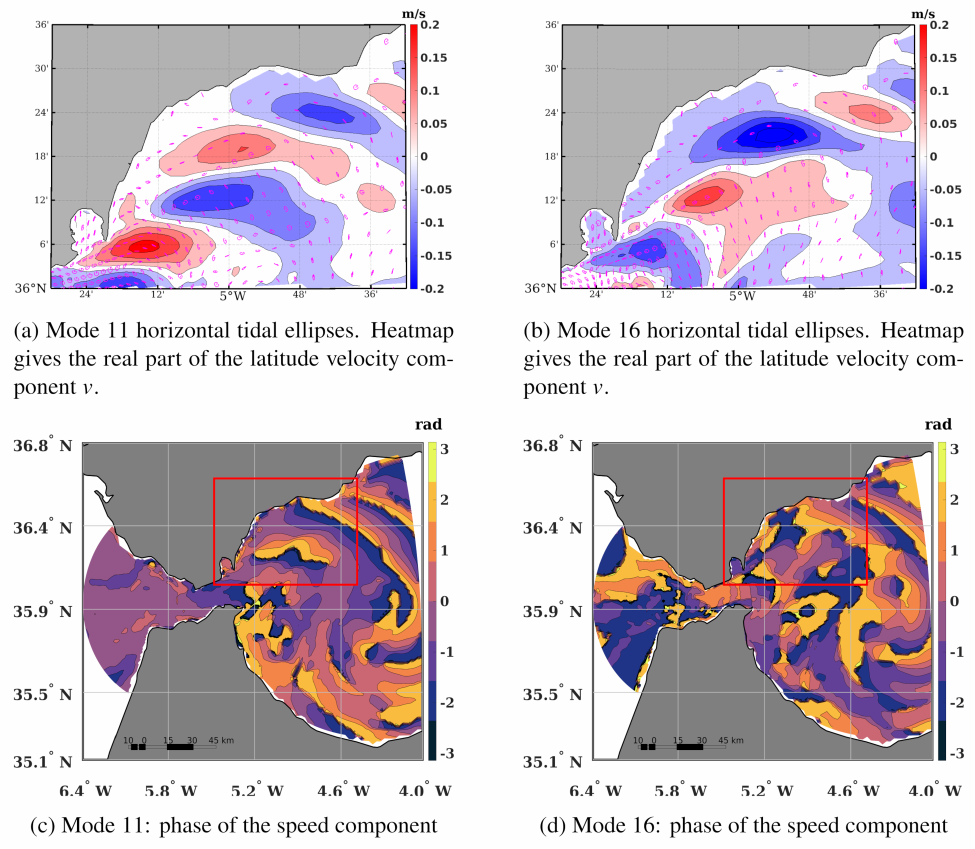}
  \caption{Modes 11 and 16. (a)-(b) Horizontal tidal ellipses of surface velocity. Heatmap shows the real part of the meridional velocity component. Arrows represents the surface velocity of mode 1. The tidal ellipses demonstrate that these modes capture meanders of the Atlantic jet, that is oscillations perpendicular to the mean direction of the jet. (c)-(d) Phase of the speed component. The parallel phase lines indicate the region where the mode propagates as a wavefront. Frames correspond to the boundaries of the plots in the top row.}
  \label{fig:R10}
\end{figure*}

Modes 11 and 16, with respective periods of 24.34h and 26.04h, are roughly-diurnal modes that decay over the duration of the simulation (Mode 16 decays at roughly twice the rate of Mode 11). While mode 11 could be harmonically related to mode 6 (as 24.34 = 2×12.12), mode 16 is unlikely a harmonic of either modes 2 or 6. These modes are most active at the eastern end of the Strait and within the Alboran Sea (Figure~\ref{fig:R9_1}). We associate Mode 11 with meanders of the Atlantic jet because its velocity and phase lines at the surface are roughly perpendicular to the direction of the mean flow along the northern coast of the Alboran Sea, as shown in Figure~\ref{fig:R10}. The surface velocity and phase associated with Mode 16 is more complex but still consistent with a meandering motion of the jet. The phase lines of each mode connect the geographical locations that simultaneously experience the peaks of the speed under each mode. They can be used to infer the speed of spatial propagation of the modal crest across the domain, as the travel time is exactly one period of oscillation between two isocurves of the same phase angle value. Parallel phase lines indicate a geographical region where the mode is propagating as a wavefront, traveling perpendicularly to the phase line. From the bottom row of Figure~\ref{fig:R10} Mode 11 captures the wavefront of the disturbance that moves across the entire length of the Atlantic jet at a speed of about 30 km/day. Mode 16, on the other hand, has a much more disorganized phase immediately to the east of the Strait; however, it also organizes as a traveling front further along the coast of Spain.

The meanders captured by modes 11 and 16 are distinguishable in velocity snapshots of the simulation, and they seem to be associated with the generation of potential vorticity on the northern coast of the strait, which is later carried by the current towards the Alboran Sea. As suggested in Figure~\ref{fig:R11}, and in a video in the Supplementary material, high potential vorticity water from the shallow, NE region of the strait, is stripped away by tidally-driven surges of surface inflow in the Mediterranean and carried with the Atlantic jet.

\begin{figure*}[htb]
  \centering
  \includegraphics[width=\textwidth]{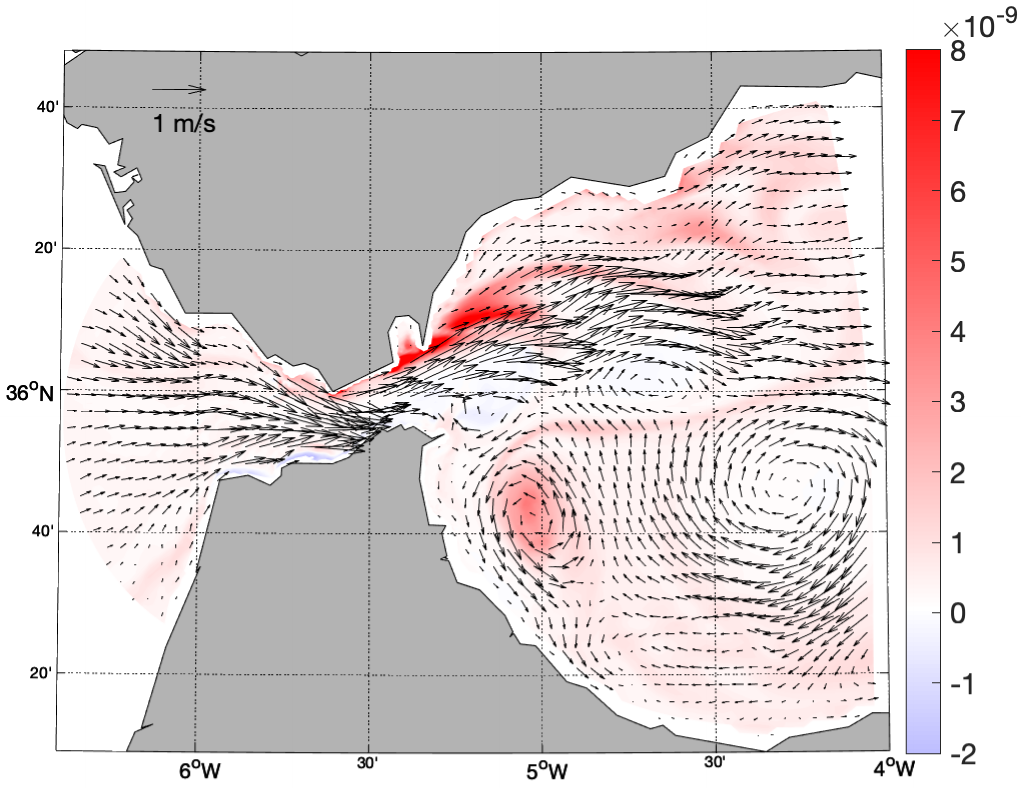}
  \caption{Snapshot of model surface velocity (arrows) and Ertel’s potential vorticity (in
PVU) at z=5m depth.}
  \label{fig:R11}
\end{figure*}

\subsection{Kelvin waves along the African coast}\label{subsec:K_AC}

\begin{figure*}[htb]
  \centering
  \includegraphics[width=\textwidth]{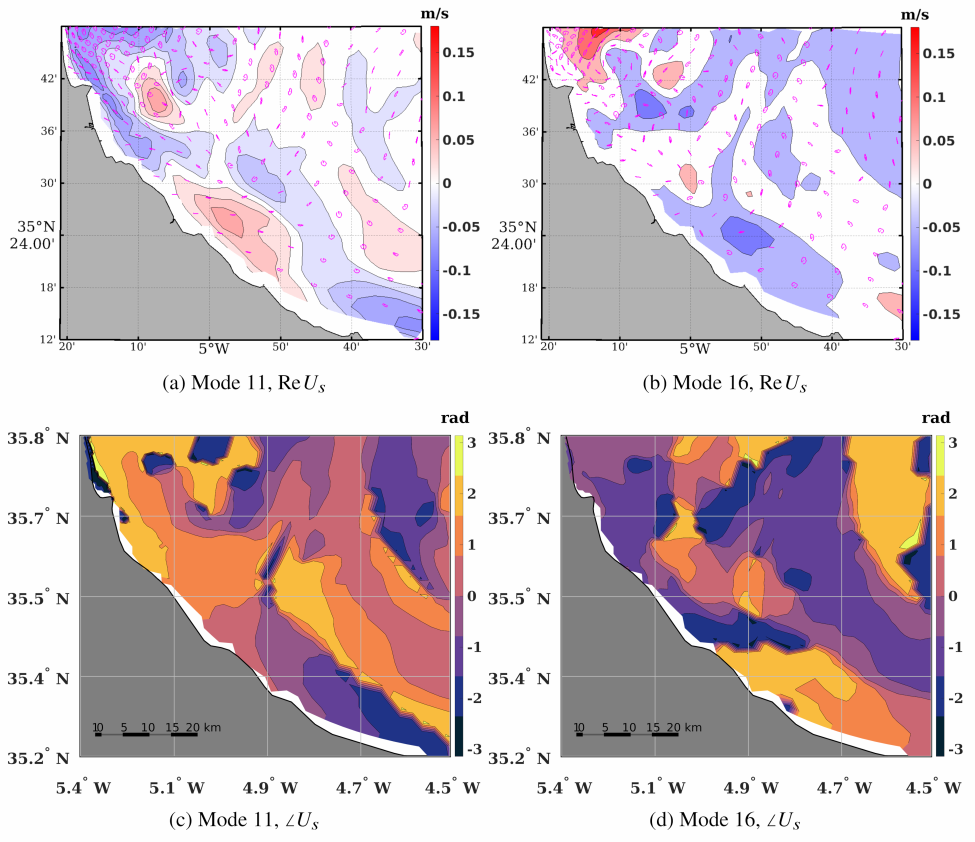}
  \caption{Near-diurnal modes 11 and 16 visualized at the surface layer, along the african coast of the Alboran sea. Top row: tidal ellipses of the horizontal velocity component; heat map represents the modes contribution to the horizontal speed. Bottom row: phase of the speed component. The eight adjacent phase lines are mapped into each other every \(\approx 24 hr/8 = 3 hr\) approximately. 
}\label{fig:R12}
\end{figure*}

Kelvin and other coastal-trapped waves can propagate along coastlines and can do so for long distances without losing their energy. In the Northern Hemisphere, these waves propagate with the coastline to the right of their propagation direction. Near-diurnal modes 11 and 16, with periods of 24.34 and 26.04 hours respectively, contain some features of coastal-trapped wave propagation. The surface tidal ellipses for mode 11 tend to be elongated near the African coast, with a general tendency to align along the coast (Figure~\ref{fig:R12}) (A Kelvin wave would have zero onshore velocity and therefore tidal ellipses that would be flat and parallel to the coast.)The phase speed of Mode 11, estimated from the phase distribution along the coast (Figure~\ref{fig:R12}) is approximately 1 m/s, which is close the theoretical estimate. 
The latter is taken as the speed of a long, first-vertical mode gravity wave speed, calculated as \( c = \sqrt{\frac{g' h_1 h_2}{h_1 + h_2}} \), where \( g' = 0.015 \, \text{m/s}^2 \) is the reduced gravity and \( h_1 = 100 \, \text{m} \), \( h_2 = 200 \, \text{m} \) are representative values of the Atlantic and Mediterranean layers thickness over the African slope, respectively, derived from the mean stratification.
In addition, the velocity magnitude of mode 11 is generally strongest within a Rossby radius of deformation, here about 30 km, of the coastline. It would thus appear that mode 11 contains at least some characteristics of a Kelvin wave propagating eastward along the coast of Africa. 

Mode 16 also exhibits a degree of coastal trapping (Figure~\ref{fig:R12}b), but here the tidal ellipses are rounder, and the phase distribution (Figure~\ref{fig:R12}d) is more complex, making it more difficult to compute an along-coast phase speed. For both modes, the offshore phase distribution is quite complex, not surprising in view of the presence of offshore topographic features. These modes were the two most robustly computed DMD modes with respect to the leave-one-out analysis, except for the time-invariant mode (Section 4~\ref{sec:robustness}). 

\subsection{Higher harmonics internal tides}\label{sec:H_HT}

\begin{figure*}[htb]
  \centering
  \includegraphics[width=\textwidth]{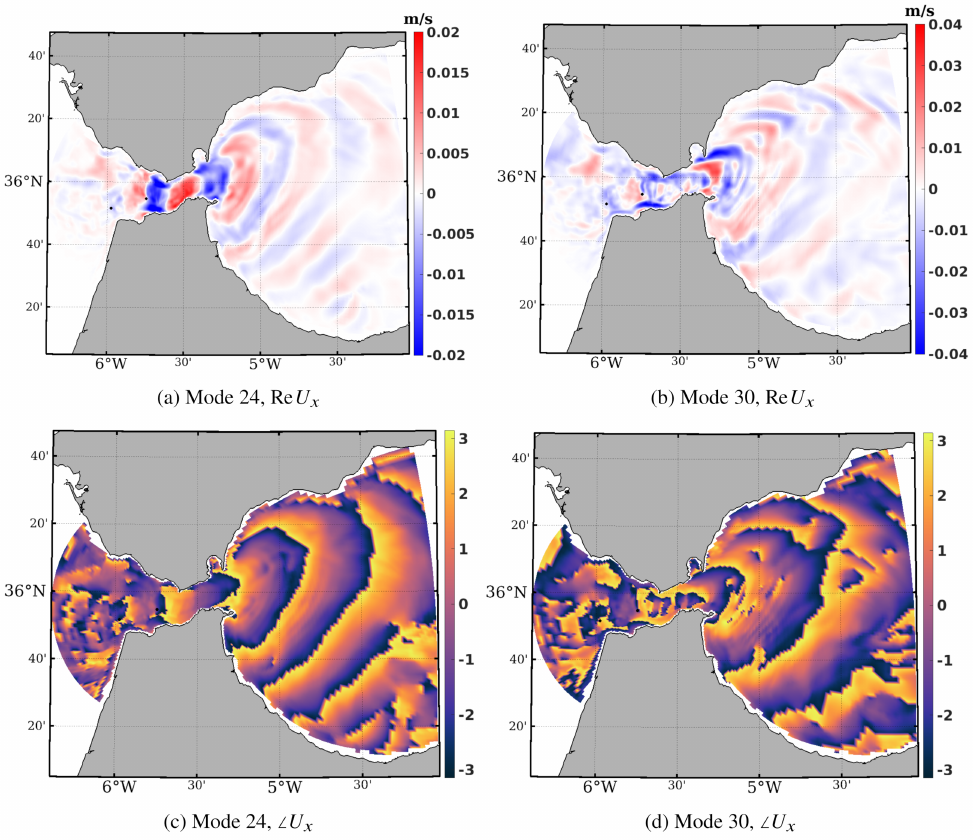}
  \caption{Modes 24 and 30 surface horizontal velocity (top) and phases (bottom).}
  \label{fig:R13}
\end{figure*}

Modes 24 and 30 are fourth-diurnal and third-diurnal, respectively, with periods of 6.19 h and 8.26 h. Mode 24 is likely the first harmonic of Mode 2 (\(6.19 \, \text{h} \approx \frac{12.34 \, \text{h}}{2}\)), while Mode 30 is the second harmonic of Mode 11 (\(8.26 \, \text{h} \approx \frac{24.34 \, \text{h}}{3}\)). Modes 24 and 30 can be thought of as nonlinear tidal constituents arising from the tidal dynamics of the Strait of Gibraltar.
The main source of nonlinearily is the CS region, where the flow undergoes a series of internal hydraulic transitions (Figure~\ref{fig:2_R2}). The nonlinear interaction between the tidal flow and the bottom topography gives rise to internal waves of harmonic frequencies that are captured by Modes 24 and 30. These modes show patterns of progressive waves in the Alboran Sea emanating from the CS (Figure~\ref{fig:R13}). The wavelength of these waves, estimated from the phase of Modes 24 and 30, is approximately 40 km, resulting in a phase velocity of $\approx 1.5 \, \text{m/s}$. This value roughly agrees with theoretical estimates of the phase velocity for long internal disturbances in the Alboran Sea (1.6 m/s; \cite{lafuente2000tide}).     

\section{Discussion and conclusions}\label{sec:discussion}

The primary objective of this study was to apply Dynamic Mode Decomposition (DMD) to analyze the complex, time-dependent oceanographic phenomena within the Strait of Gibraltar. Our findings demonstrate DMD's capability to dissect intricate fluid dynamics into physically meaningful modes, effectively capturing essential phenomena such as tidal flows, gyres,  internal and coastal-trapped waves, and meandering motions. Specifically, the analysis isolated critical dynamical features, including the meandering Atlantic Jet, tidal influences within the exchange flow at Camarinal Sill, sub-mesoscale gyres, radiating internal waves and coastal trapped waves along the African coast. Each feature is captured by distinct DMD modes with specific frequencies and decay rates. The robustness of these modes was validated through a leave-one-out analysis, assessing the sensitivity of DMD eigenvalues to the exclusion of individual data points. This analysis pinpointed clusters of eigenvalues associated with significant tidal components, underscoring the method's efficacy in isolating and analyzing key dynamic features independently of specific data points.

This study also endeavored to improve interpretability, stability, and robustness into the DMD analysis by introducing several enhancements. Key modifications include integrating horizontal speed into the snapshot vector to simplify model reduction and maintain the integrity of velocity components across dimensions. The study also reevaluates the traditional mean removal process, suggesting that omitting this step might enhance robustness for datasets with dynamic eigenvalues, supported by recent research advocating for stability benefits. Additionally, the transition from baseline least-squares to total-least-squares correction is discussed. This method employs a robust subspace estimation via the Singular Value Decomposition (SVD) of an expanded matrix, aiming to improve the accuracy of derived dynamic modes by mitigating snapshot-specific deviations.

Significant algorithmic adjustments have been implemented to enhance the stability and reliability of DMD computations. Column normalization of the snapshot matrices improves the condition number and numerical stability, essential for handling large datasets typical in oceanographic studies. The adoption of the LAPACK QR SVD algorithm over MATLAB's default SVD algorithm enhances the accuracy of computing smaller singular values and their vectors, ensuring reliable decompositions in precision-critical applications.

Rigorous robustness tests, including leave-one-out analyses, assess the stability of computed eigenvalues against minor dataset variations. These tests are crucial for ensuring the accuracy of eigenvalue calculations, which directly impact the interpretation of dynamic oceanographic phenomena such as tidal movements. By removing one snapshot at a time and recalculating the DMD, the variability in eigenvalues is observed, providing insights into the algorithm's sensitivity to data perturbations and confirming the reliability of DMD under operational conditions.

In our detailed examination of oceanographic modes, several key dynamics have been identified that play pivotal roles in shaping the oceanographic features observed in the simulation. Mode 1, characterized as an aperiodic, near-constant mode, is substantial in its alignment with the time-averaged norm of data, effectively capturing the mean properties of the Gibraltar exchange flow and those of the Western Alboran Gyre. Further insights are gained from Modes 2 and 6, both exhibiting semidiurnal periods, which dominate the unsteady dynamics at the sill of the Strait. These modes directly influence the vertical circulation patterns, such as upwelling and downwelling east of the sill, detailed in Figures~\ref{fig:R6} and \ref{fig:2_R2}. The roles of these modes are critical in mediating the the volume fluxes of water within layers, and therefore the exchange of other properties such as nutrients. 

The analysis also reveals the significant impact of Modes 11 and 16, roughly diurnal, which are intricately linked to the meandering patterns of the Atlantic jet. These modes exhibit specific phase and velocity profiles that not only enhance our understanding of jet dynamics but also suggest mechanisms for Kelvin wave propagation along the African coast, as depicted in Figures~\ref{fig:R10} and \ref{fig:R12}. Mode 11 also exhibits characteristics of Kelvin wave propagation along the coast, indicated by its elongated surface tidal ellipses, its alignment with coastal contours in Figure~\ref{fig:R12}, and the fact that is propagation speed is close to that predicted by linear wave theory. 

Modes 24 and 30, identified as higher harmonic frequencies, arise from complex nonlinear interactions at the CS region. These modes contribute to the formation of internal waves, a phenomenon critical for understanding the intricate wave patterns observed in the Alboran Sea, illustrated in Figure~\ref{fig:R13}. The generation of these internal waves through Modes 24 and 30 highlights the significant role of tidal dynamics in influencing the marine environment. 

In summary, the dynamic modes identified through our DMD analysis provide a comprehensive view of the oceanographic processes at play. From the gyre dynamics captured by Mode 1 to the intricate wave patterns facilitated by higher harmonic modes, each mode offers unique insights that are essential for a thorough understanding of the region's marine dynamics. Although these features are, of course, captured in the full numerical model, the DMD modes provided an alternative model of significantly reduced order.

\datastatement
This dataset provides 3D ocean current simulations of the Strait of Gibraltar, including zonal (U), meridional (V), and vertical (W) velocity components at grid coordinates (X, Y, Z), where X represents longitude, Y latitude, and Z depth (in meters). The simulations were conducted using the MIT General Circulation Model and consist of 144 hourly snapshots spanning six days. The data used in this study are available for access and can be obtained from \url{https://github.com/SathsaraDias/Strait_of_Gibraltar.git}. For any additional information or specific requests, please contact the authors.

\section{Funding}

Jose C. Sanchez-Garrido gratefully acknowledges funding from the regional institution Junta de Andalucía through the research project ESMER4 (PCM0005, Plan de Recuperación, Transformación y Resiliencia). Pratt was supported by NSF grant OCE\# 2124210. The research of E. Bollt is supported by the ONR, ARO, DARPA RSDN, as well as the NIH and NSF CRCNS programs.

\printunsrtglossary[type=symbols, title={\textit{List of Symbols}}]

\clearpage
\acknowledgments
Some of the computing for this project was performed on the Clarkson University ACRES cluster.
We would like to thank Clarkson University Office of Information Technology, and the National Science Foundation, who partially funded these resources under Grant No.
1925596.

%
%
\datastatement

\bibliographystyle{ametsocV6}
\bibliography{references}

\end{document}